\documentclass[12pt,twoside,a4paper,openany]{article}
\usepackage[centertags]{amsmath}
\usepackage{amsfonts}
\usepackage[all]{xy}
\usepackage{graphicx}
\usepackage{amsthm}
\usepackage{newlfont}
\usepackage[latin1]{inputenc}
\usepackage{amsmath}
\usepackage{eucal}
\usepackage{amssymb}
\usepackage{titlesec}
\titleformat{\subsection}[runin]
       {\bfseries}
			 {\thesubsection}
       {0.5em}
       {}
       []

\newtheorem{thm}{Theorem}[section]

\newtheorem{lem}[thm]{Lemma}
\newtheorem{prop}[thm]{Proposition}

\newtheorem{subs}[thm]{}
\newtheorem{rem}[thm]{Remark}

\begin{document}

\begin{center}
\textbf{\Large Limit linear series for curves of compact type}

\textbf{\Large with three irreducible components}
\end{center}

\begin{center}
\text{Gabriel Mu\~noz}
\end{center}

\begin{abstract}

Our aim in this work is to study exact Osserman limit linear series on curves of compact type $X$ with three irreducible components. This case is quite different from the case of two irreducible components studied by Osserman. For instance, for curves of compact type with two irreducible components, every refined Eisenbud-Harris limit linear series has a unique exact extension. But, for the case of three irreducible components, this property is no longer true. We find a condition characterizing when a given refined Eisenbud-Harris limit linear series has a unique exact extension. To do this, it is necessary to understand how to construct exact extensions. We find a constructive method, which describes how to construct all exact extensions of refined limit linear series. By our method, we get that every refined limit linear series has at least one exact extension.

\end{abstract}

\section{Introduction}

In Algebraic Geometry, the theory of linear series on smooth curves is closely related to that of Abel maps. The fibers of Abel maps consist precisely of complete linear series. 

For curves of compact type, (Eisenbud and Harris 1986) developed the theory of limit linear series as an analogue of linear series. This theory is very powerful for degeneration arguments on curves. The idea is to analyze how linear series degenerate when a family of smooth curves degenerates to a compact type curve. Eisenbud and Harris approached this situation by considering only the possibe limit line bundles with nonnegative multidegree and degree $d$ on one irreducible component of the curve. 
(Osserman 2006a) developed a new and more functorial construction for the theory of limit linear series. The basic idea is to consider all possible limit line bundles with nonnegative multidegree. Thus, Osserman limit linear series carry more information about limit line bundles. This new theory has a generalization to higher rank vector bundles (Osserman 2014).

Abel maps for curves of compact type have been studied by (Coelho and Pacini 2010). Recently, for curves of compact type $X$ with two irreducible components, (Esteves and Osserman 2013) related limit linear series to fibers of Abel maps via the definition of limit linear series by (Osserman 2006a). They studied the notion of exact limit linear series. These contain in particular all limits of linear series on the generic fiber in a regular smoothing family.

Also, for curves of compact type $X$ with two irreducible components, (Osserman 2006b) studied the space of limit linear series corresponding to a given Eisenbud-Harris limit linear series. He obtained an upper bound for the dimension of that space. Using this result, he also obtained a simple proof of the Brill-Noether theorem using only the limit linear series theory.
 
Our aim in this work is to study exact limit linear series on curves of compact type $X$ with three irreducible components. This case is quite different from the case of two irreducible components. For instance, for curves of compact type with two irreducible components, (Osserman 2006a) showed that every refined Eisenbud-Harris limit linear series has a unique exact extension. But, for the case of three irreducible components, this property is no longer true. 

We will study the case of exact limit linear series which are obtained as the unique exact extension of a refined Eisenbud-Harris limit linear series. For curves $X$ consisting of a chain of smooth curves, (Osserman 2014) studied the notion of chain adaptable refined limit linear series. He showed that every chain adaptable refined limit linear series has a unique extension. In particular, for our case of curves of compact type $X$ with three irreducible components, chain adaptable refined limit linear series are an example of refined Eisenbud-Harris limit linear series with a unique extension. It turns out that, by our Theorem \ref{unique}, every refined Eisenbud-Harris limit linear series having a unique extension is chain adaptable. We should mention that, for our case of curves of compact type with three irreducible components, it is easy to see that every chain adaptable refined limit linear series has a unique extension (we give a simple proof of that in Theorem \ref{unique}), but proving that the unique extension property implies the chain adaptable condition is not easy at all. We really need to understand how to construct extensions (in Section \ref{extension}, we describe a method for the construction of all exact extensions).

We now explain the contents of the paper in more detail, especially the statement of our main theorem, which is Theorem \ref{unique}. We begin with the notation of a limit linear series. In (Esteves and Osserman 2013), a limit linear series of degree $d$ and dimension $r$ on a curve of compact type $X$ with two irreducible components $Y$ and $Z$, meeting transversally at a point $P$, is a collection $\mathfrak{g}:=(\mathcal{L},V_{0},\ldots,V_{d})$, where $\mathcal{L}$ is an invertible sheaf on $X$ of degree $d$ on $Y$ and degree $0$ on $Z$, and $V_{i}$ is a vector subspace of $H^{0}(X,\mathcal{L}^{i})$ of dimension $r+1$, for each $i=0,\ldots,d$, where $\mathcal{L}^{i}$ is the invertible sheaf on $X$ with restrictions $\mathcal{L}\big|_{Y}(-iP)$ and $\mathcal{L}\big|_{Z}(iP)$, and these vector subspaces are linked by certain natural maps between the sheaves $\mathcal{L}^{i}$. Thus, a limit linear series is defined by a collection of pairs $(\mathcal{L}^{i}, V_{i})$, for each $i=0,\ldots,d$, and we can use the notation $\{(\mathcal{L}^{i}, V_{i})\}_{i}$, where $0\leq i\leq d$. For each $i=0,\ldots,d$, the invertible sheaf $\mathcal{L}^{i}$ has multidegree $\underline{d}:=(d-i,i)$. So, for each $i=0,\ldots,d$, setting $\underline{d}:=(d-i,i)$, $\mathcal{L}_{\underline{d}}=\mathcal{L}^{i}$ and $V_{\underline{d}}=V_{i}$, a limit linear series can also be denoted by a collection $\{(\mathcal{L}_{\underline{d}},V_{\underline{d}})\}_{\underline{d}}$, where $\underline{d}\geq \underline{0}$ of total degree $d$. We will use this notation for the case of three irreducible components. 

In this work, $X$ will denote the union of three smooth curves $X_1, X_2$ and $X_3$, such that $X_1$ and $X_2$ meet transversally at a point $A$, and $X_2$ and $X_3$ meet transversally at a point $B$, with $A\neq B$. A limit linear series of degree $d$ and dimension $r$ on $X$  is a collection $\mathfrak{g}:=\{(\mathcal{L}_{\underline{d}},V_{\underline{d}})\}_{\underline{d}}$, where $\underline{d}\geq \underline{0}$ of total degree $d$, each $\mathcal{L}_{\underline{d}}$ is an invertible sheaf of multidegree $\underline{d}$, and $V_{\underline{d}}$ is a vector subspace of $H^{0}(X,\mathcal{L}_{\underline{d}})$ of dimension $r+1$, for each $\underline{d}$, where $\mathcal{L}_{\underline{d}}$ is the invertible sheaf on $X$ with restrictions $\mathcal{L}\big|_{X_1}(-(d-i)A)$, $\mathcal{L}\big|_{X_2}((d-i)A-lB)$ and $\mathcal{L}\big|_{X_3}(lB)$, with $\underline{d}=(i,d-i-l,l)$ and $\mathcal{L}:=\mathcal{L}_{(d,0,0)}$, and the subspaces $V_{\underline{d}}$ are linked by certain natural maps between the sheaves $\mathcal{L}_{\underline{d}}$; see Subsection \ref{limit1} for the precise definition of limit linear series.

We restrict our attention to the case of exact limit linear series which are obtained as the unique exact extension of a refined Eisenbud-Harris limit linear series. Given a refined limit linear series, for any $i\geq 0$ and $l\geq 0$ such that $i+l\leq d$ there is a natural space, here denoted $K_{il}$, satisfying that, for every extension $\mathfrak{g}:=\{(\mathcal{L}_{\underline{d}},V_{\underline{d}})\}_{\underline{d}}$, $K_{il}$ contains $V_{(i,d-i-l,l)}$. It turns out that these spaces are very important to understand when the extension is unique (spaces with the same property as the $K_{il}$ appeared in (Osserman 2006a) and (Osserman 2014), for curves of compact type with two irreducible components and for curves of compact type with more than two components, respectively).

We recall the analogous spaces for the case of curves of compact type $X$ with two components. Given an invertible sheaf $\mathcal{L}$ on $X$ of degree $d$ on $Y$ and degree $0$ on $Z$, for each nonnegative multidegree $\underline{d}:=(i,d-i)$, let $\mathcal{L}_{\underline{d}}$ be the invertible sheaf on $X$ with restrictions $\mathcal{L}\big|_{Y}(-(d-i)P)$ and $\mathcal{L}\big|_{Z}((d-i)P)$, and let $V_{Y}$ and $V_{Z}$ be vector subspaces of $H^{0}(\mathcal{L}_{(d,0)})$ and $H^{0}(\mathcal{L}_{(0,d)})$, respectively, of dimension $r+1$, such that $\{(\mathcal{L}_{(d,0)}\big|_{Y},V_{Y}\big|_{Y}), (\mathcal{L}_{(0,d)}\big|_{Z}, V_{Z}\big|_{Z})\}$ is a refined Eisenbud-Harris limit linear series. In this situation, (Osserman 2006a) showed that the unique exact extension is given as follows: 

$(*)$ For each nonnegative multidegree $\underline{d}:=(i,d-i)$, $V_{\underline{d}}$ is the space of sections $s$ of $H^{0}(X,\mathcal{L}_{\underline{d}})$ whose image in $H^{0}(X,\mathcal{L}_{(d,0)})$ belongs to $V_{Y}$ and vanishes at $P$ with order at least $d-i$, and whose image in $H^{0}(X,\mathcal{L}_{(0,d)})$ belongs to $V_{Z}$ and vanishes at $P$ with order at least $i$.

Now, in our case of three irreducible components, consider the analogous situation: given an invertible sheaf $\mathcal{L}$ on $X$ of multidegree $(d,0,0)$, for each nonnegative multidegree $\underline{d}:=(i,d-i-l,l)$, let $\mathcal{L}_{\underline{d}}$ be the invertible sheaf on $X$ with restrictions $\mathcal{L}\big|_{X_{1}}(-(d-i)A)$, $\mathcal{L}\big|_{X_{2}}((d-i)A-lB)$ and $\mathcal{L}\big|_{X_{3}}(lB)$, and let $V_{X_1}$, $V_{X_2}$ and $V_{X_3}$ be vector subspaces of $H^{0}(\mathcal{L}_{(d,0,0)})$, $H^{0}(\mathcal{L}_{(0,d,0)})$ and $H^{0}(\mathcal{L}_{(0,0,d)})$, respectively, of dimension $r+1$, such that $\{(\mathcal{L}_{(d,0,0)}\big|_{X_{1}},V_{X_1}\big|_{X_{1}}), (\mathcal{L}_{(0,d,0)}\big|_{X_{2}}, V_{X_2}\big|_{X_{2}}), (\mathcal{L}_{(0,0,d)}\big|_{X_{3}},V_{X_3}\big|_{X_{3}})\}$ is a refined Eisenbud-Harris limit linear series. For each $i,l$, define $K_{il}$ as the natural generalization of $(*)$. More specifically:

For each nonnegative multidegree $\underline{d}:=(i,d-i-l,l)$, $K_{il}$ is the space of sections $s$ of $H^{0}(X,\mathcal{L}_{\underline{d}})$ whose image in $H^{0}(X,\mathcal{L}_{(d,0,0)})$ belongs to $V_{X_1}$ and vanishes at $A$ with order at least $d-i$, whose image in $H^{0}(X,\mathcal{L}_{(0,d,0)})$ belongs to $V_{X_2}$ and vanishes at $A$ with order at least $i$ and vanishes at $B$ with order at least $l$, and whose image in $H^{0}(X,\mathcal{L}_{(0,0,d)})$ belongs to $V_{X_3}$ and vanishes at $B$ with order at least $d-l$; see Subsection \ref{kernels1} for the precise definition of the spaces $K_{il}$.

The difference with the case of two irreducible components is the fact that in our case the spaces defined above does not necessarily have dimension $r+1$. 

In Section \ref{extension}, we describe a method for the construction of all exact extensions of refined limit linear series. As a consequence, we get that every refined limit linear series has at least one exact extension. We use the method of Section \ref{extension} to understand when a refined limit linear series has a unique exact extension. Our Theorem \ref{unique} says that, for a refined limit linear series:

{\em There is a unique exact extension if and only if }
\begin{center}
\text{dim}\,$K_{il}=r+1$ $if$ $i+l\leq d$, $b_{j-1}<i\leq b_{j}$, $b'_{k-1}<l\leq b'_{k}$ $and$ $j+k\leq r+1$,
\end{center}
{\em where $b_0,\ldots,b_r$ are the orders of vanishing at $A$ and $b'_0,\ldots,b'_r$ are the orders of vanishing at $B$, all orders correspond to the linear series on $X_{2}$.
Moreover, this condition is also equivalent to the existence of a unique extension.}

It follows from the exact sequence defining $K_{il}$ (see Subsection \ref{kernels1}), that the condition in Theorem \ref{unique} is equivalent to the following condition:
 \begin{center}
\text{dim}\,$V_{X_2}(-iA-lB)=r+1-j-k$
\end{center}
if $i+l\leq d$, $b_{j-1}<i\leq b_{j}$, $b'_{k-1}<l\leq b'_{k}$ and $j+k\leq r+1$.
This condition means exactly that our refined limit linear series is chain adaptable. 

For compact type curves with two irreducible components, given an exact limit linear series, (Esteves and Osserman 2013) associated a closed subscheme of the fiber of the corresponding Abel map. In a subsequent work, for compact type curves with three irreducible components, we will give a description of this closed subscheme when the underlying exact limit linear series is the unique extension of a refined limit linear series.
 
Our techniques can be generalized to the case of compact type curves with an arbitrary number of irreducible components (work in progress).

\section{Preliminaries}\label{limit}
\begin{subs}{\em(}Limit linear series{\em)}\label{limit1} {\em Throughout this article, $X$ will denote the union of three smooth curves $X_1, X_2$ and $X_3$, such that $X_1$ and $X_2$ meet transversally at a point $A$, and $X_2$ and $X_3$ meet transversally at a point $B$, with $A\neq B$. If $Y$ is a reduced union of some components of $X$, we get the following exact sequence
\[0\rightarrow \mathcal{L}\big|_{Y^{c}}(-Y\cap Y^{c})\rightarrow \mathcal{L}\rightarrow \mathcal{L}\big|_{Y}\rightarrow 0,\]
for any invertible sheaf $\mathcal{L}$ on $X$.
Let $\mathcal{L}$ be an invertible sheaf on $X$ of degree $d$ on $X_1$ and degree $0$ on $X_2$ and $X_3$. For any $i\geq 0$ and $l\geq 0$ such that $i+l\leq d$, let $\mathcal{L}_{(i,d-i-l,l)}$ be the invertibe sheaf on $X$ with restrictions $\mathcal{L}\big|_{X_1}(-(d-i)A)$, $\mathcal{L}\big|_{X_2}((d-i)A-lB)$ and $\mathcal{L}\big|_{X_3}(lB)$. Note that $\mathcal{L}_{(i,d-i-l,l)}$ has multidegree $(i,d-i-l,l)$. For any $i\geq 0$ and $l\geq 0$ such that $i+l\leq d$, let $\underline{d}:=(i,d-i-l,l)$ and set
 \[\widetilde{\underline{d}}:=\left\{
            \begin{array}{rcl} (i-1,d-i-l+1,l) & \mbox{if} & q=1,\\
               (i+1,d-i-l-2,l+1) & \mbox{if} & q=2,\\
							 (i,d-i-l+1,l-1) & \mbox{if} & q=3.
                           \end{array}\right.\]


Whenever $\widetilde{\underline{d}}\geq \underline{0}$, there are natural maps
\[\varphi_{\underline{d},\widetilde{\underline{d}}}: \mathcal{L}_{\underline{d}}\rightarrow \mathcal{L}_{\underline{d}}\big|_{X_{q}^{c}}=\mathcal{L}_{\widetilde{\underline{d}}}\big|_{X_{q}^{c}}(-X_{q}\cap X_{q}^{c})\hookrightarrow \mathcal{L}_{\widetilde{\underline{d}}},\]
\[\varphi_{\widetilde{\underline{d}},\underline{d}}: \mathcal{L}_{\widetilde{\underline{d}}}\rightarrow \mathcal{L}_{\widetilde{\underline{d}}}\big|_{X_{q}}=\mathcal{L}_{\underline{d}}\big|_{X_{q}}(-X_{q}\cap X_{q}^{c})\hookrightarrow \mathcal{L}_{\underline{d}},\]
where the first map in each composition is the restriction map and the last maps are the natural inclusions. Note that the compositions $\varphi_{\underline{d},\widetilde{\underline{d}}} \varphi_{\widetilde{\underline{d}},\underline{d}}$ and $\varphi_{\widetilde{\underline{d}},\underline{d}} \varphi_{\underline{d},\widetilde{\underline{d}}}$ are zero.

If $Y$ is a subcurve of $X$, for any $\underline{d}$ and for any subspace $V\subseteq H^{0}(X,\mathcal{L}_{\underline{d}})$, we denote by $V^{Y,0}$ the subspace of $V$ of sections that vanish on $Y$. If $Y$ is an irreducible component of $X$, we denote by $\mathcal{L}_{Y}$ the invertible sheaf $\mathcal{L}_{\underline{d}}$, where the component of $\underline{d}$ corresponding to $Y$ is equal to $d$ and the other components of $\underline{d}$ are $0$. Also, to ease notation, let $\mathcal{L}_{il}:=\mathcal{L}_{(i,d-i-l,l)}$.

Fix integers $d$ and $r$. A {\em limit linear series} on $X$ of degree $d$ an dimension $r$ is a collection consisting of an invertible sheaf $\mathcal{L}$ on $X$ of degree $d$ on $X_1$ and degree $0$ on $X_2$ and $X_3$, and vector subspaces $V_{\underline{d}}\subseteq H^{0}(X,\mathcal{L}_{\underline{d}})$ of dimension $r+1$, for each $\underline{d}:=(i,d-i-l,l)\geq \underline{0}$, such that $\varphi_{\underline{d},\widetilde{\underline{d}}}(V_{\underline{d}})\subseteq V_{\widetilde{\underline{d}}}$ and $\varphi_{\widetilde{\underline{d}},\underline{d}}(V_{\widetilde{\underline{d}}})\subseteq V_{\underline{d}}$ for each $\underline{d}\geq \underline{0}$, whenever $\widetilde{\underline{d}}\geq \underline{0}$.

Given a limit linear series, if $Y$ is an irreducible component of $X$, we denote by $V_{Y}$ the corresponding subspace of $H^{0}(X,\mathcal{L}_{Y})$. Also, we denote by $V_{il}$ the corresponding subspace of $H^{0}(X,\mathcal{L}_{il})$.

The conditions $\varphi_{\underline{d},\widetilde{\underline{d}}}(V_{\underline{d}})\subseteq V_{\widetilde{\underline{d}}}$ and $\varphi_{\widetilde{\underline{d}},\underline{d}}(V_{\widetilde{\underline{d}}})\subseteq V_{\underline{d}}$ are called {\em the linking condition}, and we say that $V_{\underline{d}}$ and $V_{\widetilde{\underline{d}}}$ are linked by the maps $\varphi_{\underline{d},\widetilde{\underline{d}}}$ and $\varphi_{\widetilde{\underline{d}},\underline{d}}$. 

Note that $\varphi_{\underline{d},\widetilde{\underline{d}}}: V_{\underline{d}}\rightarrow V_{\widetilde{\underline{d}}}$ has kernel $V_{\underline{d}}^{X_{q}^{c},0}$ and image contained in $V_{\widetilde{\underline{d}}}^{X_{q},0}$. Analogously, the map $\varphi_{\widetilde{\underline{d}},\underline{d}}: V_{\widetilde{\underline{d}}}\rightarrow V_{\underline{d}}$ has kernel $V_{\widetilde{\underline{d}}}^{X_{q},0}$ and image contained in $V_{\underline{d}}^{X_{q}^{c},0}$.

A limit linear series $\{(\mathcal{L}_{\underline{d}},V_{\underline{d}})\}_{\underline{d}}$ is called {\em exact} if
\begin{center}
Im\,$(\varphi_{\underline{d},\widetilde{\underline{d}}}: V_{\underline{d}}\rightarrow V_{\widetilde{\underline{d}}})=V_{\widetilde{\underline{d}}}^{X_{q},0}$ and 
Im\,$(\varphi_{\widetilde{\underline{d}},\underline{d}}: V_{\widetilde{\underline{d}}}\rightarrow V_{\underline{d}})=V_{\underline{d}}^{X_{q}^{c},0}$
\end{center}
for each $\underline{d}:=(i,d-i-l,l)\geq \underline{0}$, whenever $\widetilde{\underline{d}}\geq \underline{0}$.

}
\end{subs}

\begin{rem}\label{anulamiento}
{\em Note that, from the construction of the invertible sheaves $\mathcal{L}_{il}$ and the maps $\varphi_{\underline{d},\widetilde{\underline{d}}}$, we have that, for $\widetilde{q}\neq q$, $s\in H^{0}(\mathcal{L}_{\underline{d}})$ vanishes on $X_{\widetilde{q}}$ if and only if $\varphi_{\underline{d},\widetilde{\underline{d}}}(s)\in H^{0}(\mathcal{L}_{\widetilde{\underline{d}}})$ vanishes on $X_{\widetilde{q}}$. In particular, given a limit linear series $\{(\mathcal{L}_{\underline{d}},V_{\underline{d}})\}_{\underline{d}}$, we get natural inclusions, for $\widetilde{q}\neq q$
\[V_{\underline{d}}/V_{\underline{d}}^{X_{\widetilde{q}},0}\hookrightarrow V_{\widetilde{\underline{d}}}/V_{\widetilde{\underline{d}}}^{X_{\widetilde{q}},0}.\]
Also, if $s\in H^{0}(\mathcal{L}_{\widetilde{\underline{d}}})$ and $\varphi_{\widetilde{\underline{d}},\underline{d}}(s)\in H^{0}(\mathcal{L}_{\underline{d}})$ vanishes on $X_{q}$, then $\varphi_{\widetilde{\underline{d}},\underline{d}}(s)$ vanishes on $X_{q}\cup X_{q}^{c}=X$, as $\varphi_{\widetilde{\underline{d}},\underline{d}}(H^{0}(\mathcal{L}_{\widetilde{\underline{d}}}))$ is contained in the kernel of the map $H^{0}(\mathcal{L}_{\underline{d}})\rightarrow H^{0}(\mathcal{L}_{\underline{d}}\big|_{X_{q}^{c}})$, and hence $\varphi_{\widetilde{\underline{d}},\underline{d}}(s)=0$, which implies that $s$ vanishes on $X_{q}$. Thus, given a limit linear series $\{(\mathcal{L}_{\underline{d}},V_{\underline{d}})\}_{\underline{d}}$, we get the natural inclusion
\[V_{\widetilde{\underline{d}}}/V_{\widetilde{\underline{d}}}^{X_{q},0}\hookrightarrow V_{\underline{d}}/V_{\underline{d}}^{X_{q},0}.\]
}
\end{rem}

\section{The kernel $K_{il}$}\label{kernels}
\begin{subs}{\em(}The kernel $K_{il}${\em)}\label{kernels1} {\em Let $\mathcal{L}$ be an invertible sheaf on $X$ of degree $d$ on $X_1$ and degree $0$ on $X_2$ and $X_3$. For any $i\geq 0$ and $l\geq 0$ such that $i+l\leq d$, recall that $\mathcal{L}_{il}$ denotes the invertibe sheaf on $X$ with restrictions $\mathcal{L}\big|_{X_1}(-(d-i)A)$, $\mathcal{L}\big|_{X_2}((d-i)A-lB)$ and $\mathcal{L}\big|_{X_3}(lB)$. 

Note that
\begin{equation*}
\mathcal{L}_{il}\big|_{X_1}=\mathcal{L}_{X_{1}}\big|_{X_1}(-(d-i)A),\\ 
\mathcal{L}_{il}\big|_{X_2}=\mathcal{L}_{X_{2}}\big|_{X_2}(-iA-lB), \\
\mathcal{L}_{il}\big|_{X_3}=\mathcal{L}_{X_{3}}\big|_{X_3}(-(d-l)B).
\end{equation*}
Then, we get the natural exact sequence:

\[0\rightarrow H^{0}(\mathcal{L}_{il})\rightarrow H^{0}(\mathcal{L}_{X_{1}}\big|_{X_1}(-(d-i)A))\oplus H^{0}(\mathcal{L}_{X_2}\big|_{X_2}(-iA-lB))\]
\[\oplus H^{0}(\mathcal{L}_{X_3}\big|_{X_3}(-(d-l)B))\rightarrow k\oplus k,\]
where the first summand in $k\oplus k$ corresponds to the point $A$ and the second summand corresponds to the point $B$. The last map in the exact sequence will be denoted $ev^{il}$.

Let $V_{X_1}, V_{X_2}, V_{X_3}$ be $r+1$-dimensional subspaces of $H^{0}(\mathcal{L}_{X_1}), H^{0}(\mathcal{L}_{X_2})$ and $H^{0}(\mathcal{L}_{X_3})$, respectively, such that they satisfy the linking condition. Assume that the associated Eisenbud-Harris limit linear series on $X$ is refined. Call $\mathfrak{h}$ this limit linear series.
For any $i\geq 0$ and $l\geq 0$ such that $i+l\leq d$, we define $K_{il}$ by the exact sequence:
\[0\rightarrow K_{il}\rightarrow V_{X_1}(-(d-i)A)\oplus V_{X_2}(-iA-lB)\oplus V_{X_3}(-(d-l)B)\rightarrow k\oplus k\]
Thus
\[K_{il}=(\alpha_{1 \underline{d}})^{-1}(V_{X_1}(-(d-i)A))\cap (\alpha_{2 \underline{d}})^{-1}(V_{X_2}(-iA-lB))\cap (\alpha_{3 \underline{d}})^{-1}(V_{X_3}(-(d-l)B))\]
where $\underline{d}:=(i,d-i-l,l)$, and the natural map
\begin{center}
$\alpha_{q \underline{d}}: H^{0}(\mathcal{L}_{\underline{d}})\rightarrow H^{0}(\mathcal{L}_{X_{q}})$ is the restriction to $X_{q}$,
\end{center}
for each $q=1,2,3$.
Denote by $b_0,\ldots,b_r$ the orders of vanishing of $V_{X_2}$ at $A$, and denote by $b'_0,\ldots,b'_r$ the orders of vanishing of $V_{X_2}$ at $B$. Also, let $a_0,\ldots,a_r$ denote the orders of vanishing of $V_{X_1}$ at $A$, and $c_0,\ldots,c_r$ the orders of vanishing of $V_{X_3}$ at $B$. Throughout this article, the data of this subsection will remain fixed.
}
\end{subs}
\begin{rem}\label{ultrarefinado determina los V}

{\em Note that, if $\{(\mathcal{L}_{\underline{d}},V_{\underline{d}})\}_{\underline{d}}$ is a limit linear series which is an extension of $\mathfrak{h}$, then $V_{il}\subseteq K_{il}$ for any $i\geq 0$ and $l\geq 0$ such that $i+l\leq d$.

Indeed, let $\underline{d}:=(i,d-i-l,l)$. By the linking condition, we have that $\alpha_{1 \underline{d}}(V_{\underline{d}})\subseteq V_{X_{1}}$. Since Im$(\alpha_{1 \underline{d}})\subseteq H^{0}(\mathcal{L}_{X_1}\big|_{X_1}(-(d-i)A))$, $\alpha_{1 \underline{d}}(V_{\underline{d}})\subseteq H^{0}(\mathcal{L}_{X_1}\big|_{X_1}(-(d-i)A))$, and hence $\alpha_{1 \underline{d}}(V_{\underline{d}})\subseteq V_{X_1}(-(d-i)A)$. Thus $V_{\underline{d}}\subseteq (\alpha_{1 \underline{d}})^{-1}(V_{X_1}(-(d-i)A))$. Analogously, we have $V_{\underline{d}}\subseteq (\alpha_{2 \underline{d}})^{-1}(V_{X_2}(-iA-lB))$ and $V_{\underline{d}}\subseteq (\alpha_{3 \underline{d}})^{-1}(V_{X_3}(-(d-l)B))$. It follows that $V_{\underline{d}}\subseteq K_{il}$.
}
\end{rem}

\begin{rem}\label{imagenev}
{\em By abuse of notation, denote the restriction of $ev^{il}$ to the vector subspace 
\begin{center}
$V_{X_1}(-(d-i)A)\oplus V_{X_2}(-iA-lB)\oplus V_{X_3}(-(d-l)B)$
\end{center}
by $ev^{il}$ as well. Notice that
\begin{equation*}
\begin{split}
(1)\,\,\,\, & (1,0)\in \text{Im}(ev^{il})\,\,\, \text{if}\,\,\, V_{X_1}(-(d-i)A)\neq V_{X_1}(-(d-i+1)A)\\
& \text{or}\,\,\, V_{X_2}(-iA-lB)\neq V_{X_2}(-(i+1)A-lB),
\end{split}
\end{equation*}
and
\begin{equation*}
\begin{split}
(2)\,\,\,\, & \text{Im}(ev^{il})\subseteq \{0\}\oplus k \,\,\,\text{if}\,\,\, V_{X_1}(-(d-i)A)=V_{X_1}(-(d-i+1)A)\\ 
& \text{and}\,\,\, V_{X_2}(-iA-lB)=V_{X_2}(-(i+1)A-lB).
\end{split}
\end{equation*}
Analogously, we have that 
\begin{equation*}
\begin{split}
(i)\,\,\,\, & (0,1)\in \text{Im}(ev^{il})\,\,\, \text{if}\,\,\, V_{X_3}(-(d-l)B)\neq V_{X_3}(-(d-l+1)B)\\
& \text{or}\,\,\, V_{X_2}(-iA-lB)\neq V_{X_2}(-iA-(l+1)B),
\end{split}
\end{equation*}
and
\begin{equation*}
\begin{split}
(ii)\,\,\,\, & \text{Im}(ev^{il})\subseteq k\oplus \{0\} \,\,\,\text{if}\,\,\, V_{X_3}(-(d-l)B)=V_{X_3}(-(d-l+1)B)\\ 
& \text{and}\,\,\, V_{X_2}(-iA-lB)=V_{X_2}(-iA-(l+1)B).
\end{split}
\end{equation*}


}
\end{rem}

\begin{rem}\label{sumaordenes}
{\em Let $C$ be a smooth curve, $L$ an invertible sheaf on $C$ of degree $d$, and $V\subseteq H^{0}(L)$ a linear series. Let $r+1:=$dim\,$V$, and let $Q_1,Q_2\in C$ distinct points. Let $e_1,\ldots,e_r$ be the orders of vanishing of $V$ at $Q_1$, and $e'_0,\ldots,e'_r$ the orders of vanishing of $V$ at $Q_2$. Then $e_j+e'_k\leq d$ if $j+k\leq r$. Furthermore, $\text{dim}\,V(-e_{j}Q_{1}-e'_{k}Q_{2})\geq r+1-(j+k)$ for any $j,k$.

Indeed,
\begin{equation*}
\begin{split}
\text{dim}\,V(-e_{j}Q_{1}-e'_{k}Q_{2})&=\text{dim}\,V(-e_{j}Q_{1})+\text{dim}\,V(e'_{k}Q_{2})-\text{dim}\,(V(-e_{j}Q_{1})+V(e'_{k}Q_{2}))\\
&\geq \text{dim}\,V(-e_{j}Q_{1})+\text{dim}\,V(e'_{k}Q_{2})-(r+1)\\
&=(r+1-j)+(r+1-k)-(r+1)=r+1-(j+k).
\end{split}
\end{equation*}
Thus, if $j+k\leq r$, then $\text{dim}\,V(-e_{j}Q_{1}-e'_{k}Q_{2})\geq 1$, which implies $h^{0}(L(-e_{j}Q_{1}-e'_{k}Q_{2}))\geq 1$, and hence deg$(L(-e_{j}Q_{1}-e'_{k}Q_{2}))\geq 0$, i.e. $e_j+e'_k\leq d$.
}
\end{rem}

\begin{prop}\label{ultrarefinado es lls}
The following statements hold:

1. {\em dim\,}$K_{il}\geq r+1$. Furthermore, {\em dim\,}$K_{il}=r+1$ if $i\leq b_0$ or $l\leq b'_0$.

2. The subspaces $K_{il}\subseteq H^{0}(\mathcal{L}_{il})$ satisfy the linking condition. 
\end{prop}
{\em Proof.} We will first prove that dim\,$K_{il}\geq r+1$. There are five cases to consider.\\
{\em Case 1:} If $i=b_j$ and $l=b'_k$ for some $j,k$.

Consider the exact sequence
\[0\rightarrow K_{il}\rightarrow V_{X_1}(-(d-i)A)\oplus V_{X_2}(-iA-lB)\oplus V_{X_3}(-(d-l)B)\rightarrow k\oplus k\]
Since
\begin{center}
$V_{X_1}(-(d-i)A)\neq V_{X_1}(-(d-i+1)A)$ and $V_{X_3}(-(d-l)B)\neq V_{X_3}(-(d-l+1)B)$,
\end{center}
as $d-i=d-b_j=a_{r-j}$ is an order of vanishing of $V_{X_1}$ at $A$ and $d-l=d-b'_k=c_{r-k}$ is an order of vanishing of $V_{X_3}$ at $B$, we have that $\text{Im}(ev^{il})=k\oplus k$, by Remark \ref{imagenev}, and hence
\begin{equation*}
\begin{split}
\text{dim}\,K_{il}&=\text{dim}\,V_{X_1}(-(d-i)A)+\text{dim}\,V_{X_2}(-iA-lB)+\text{dim}\,V_{X_3}(-(d-l)B)-2\\
&=(r+1-(r-j))+\text{dim}\,V_{X_2}(-iA-lB)+(r+1-(r-k))-2\\
&=j+k+\text{dim}\,V_{X_2}(-iA-lB)\geq j+k+(r+1-(j+k))=r+1,
\end{split}
\end{equation*}
where in the last inequality we used Remark \ref{sumaordenes}.\\
{\em Case 2:} If $b_{j-1}<i<b_{j}$ and $l=b'_k$ for some $j,k$.

Since
\begin{center}
$V_{X_1}(-(d-i)A)=V_{X_1}(-(d-i+1)A)$, $V_{X_2}(-iA-lB)=V_{X_2}(-(i+1)A-lB)$\\
and $V_{X_3}(-(d-l)B)\neq V_{X_3}(-(d-l+1)B)$,
\end{center}
as $a_{r-j}<d-i<a_{r+1-j}$ is not an order of vanishing of $V_{X_1}$ at $A$, $b_{j-1}<i<b_{j}$ is not an order of vanishing of $V_{X_2}$ at $A$ and $d-l=c_{r-k}$ is an order of vanishing of $V_{X_3}$ at $B$, we have that $\text{Im}(ev^{il})=\{0\}\oplus k$, by Remark \ref{imagenev}, and hence
\begin{equation*}
\begin{split}
\text{dim}\,K_{il}&=\text{dim}\,V_{X_1}(-(d-i)A)+\text{dim}\,V_{X_2}(-iA-lB)+\text{dim}\,V_{X_3}(-(d-l)B)-1\\
&=(r+1-(r+1-j))+\text{dim}\,V_{X_2}(-iA-lB)+(r+1-(r-k))-1\\
&=j+k+\text{dim}\,V_{X_2}(-iA-lB)\\
&=j+k+\text{dim}\,V_{X_2}(-b_{j}A-lB)\geq j+k+(r+1-(j+k))=r+1,
\end{split}
\end{equation*}
where in the last equality we used that $V_{X_2}(-iA)=V_{X_2}(-b_{j}A)$, as $b_{j-1}<i<b_{j}$, and in the last inequality we used Remark \ref{sumaordenes}.\\
{\em Case 3:} If $i=b_j$ and $b'_{k-1}<l<b'_k$ for some $j,k$.

This case is analogous to Case 2.\\
{\em Case 4:} If $b_{j-1}<i<b_{j}$ and $b'_{k-1}<l<b'_k$ for some $j,k$.

Since
\begin{center}
$V_{X_1}(-(d-i)A)=V_{X_1}(-(d-i+1)A)$, $V_{X_2}(-iA-lB)=V_{X_2}(-(i+1)A-lB)$,\\
$V_{X_2}(-iA-lB)=V_{X_2}(-iA-(l+1)B)$ and $V_{X_3}(-(d-l)B)=V_{X_3}(-(d-l+1)B)$,
\end{center}
as $d-i$ is not an order of vanishing of $V_{X_1}$ at $A$, $i$ is not an order of vanishing of $V_{X_2}$ at $A$, $l$ is not an order of vanishing of $V_{X_2}$ at $B$ and $d-l$ is not an order of vanishing of $V_{X_3}$ at $B$, we have that $\text{Im}(ev^{il})=\{0\}\oplus \{0\}$, by Remark \ref{imagenev}, and hence
\begin{equation*}
\begin{split}
\text{dim}\,K_{il}&=\text{dim}\,V_{X_1}(-(d-i)A)+\text{dim}\,V_{X_2}(-iA-lB)+\text{dim}\,V_{X_3}(-(d-l)B)\\
&=j+k+\text{dim}\,V_{X_2}(-iA-lB)\\
&=j+k+\text{dim}\,V_{X_2}(-b_{j}A-b'_{k}B)\geq j+k+(r+1-(j+k))=r+1,
\end{split}
\end{equation*}
where in the last inequality we used Remark \ref{sumaordenes}, and in the last equality we used that $V_{X_2}(-iA)=V_{X_2}(-b_{j}A)$ and $V_{X_2}(-lB)=V_{X_2}(-b'_{k}B)$.\\
{\em Case 5:} If $i<b_0$ or $i>b_r$ or $l<b'_0$ or $l>b'_r$.

We will only prove the stated inequality in the case $i<b_0$, as the other cases are analogous. Suppose $i<b_0$. Then $d-i>a_r$, and hence $V_{X_!}(-(d-i)A)=0$. Also, we have $V_{X_2}(-iA-lB)=V_{X_2}(-(i+1)A-lB)$, as $i$ is not an order of vanishing of $V_{X_2}$ at $A$. 

Suppose first that $l=b'_k$ for some $k$. Then $V_{X_3}(-(d-l)B)\neq V_{X_3}(-(d-l+1)B)$, and it follows from the exact sequence defining $K_{il}$ that
\begin{equation*}
\begin{split}
\text{dim}\,K_{il}&=\text{dim}\,V_{X_1}(-(d-i)A)+\text{dim}\,V_{X_2}(-iA-lB)+\text{dim}\,V_{X_3}(-(d-l)B)-1\\
&=0+\text{dim}\,V_{X_2}(-iA-lB)+(r+1-(r-k))-1\\
&=k+\text{dim}\,V_{X_2}(-iA-lB)\\
&=k+\text{dim}\,V_{X_2}(-lB)=k+(r+1-k)=r+1.
\end{split}
\end{equation*}
We used above that $V_{X_2}(-iA)=V_{X_2}$, as $i<b_0$.

Suppose $b'_{k-1}<l<b'_k$ for some $k$. Then $V_{X_3}(-(d-l)B)=V_{X_3}(-(d-l+1)B)$ and $V_{X_2}(-iA-lB)=V_{X_2}(-iA-(l+1)B)$, and hence
\begin{equation*}
\begin{split}
\text{dim}\,K_{il}&=\text{dim}\,V_{X_1}(-(d-i)A)+\text{dim}\,V_{X_2}(-iA-lB)+\text{dim}\,V_{X_3}(-(d-l)B)\\
&=k+\text{dim}\,V_{X_2}(-iA-lB)\\
&=k+\text{dim}\,V_{X_2}(-lB)=k+(r+1-k)=r+1.
\end{split}
\end{equation*}

Now, suppose $l<b'_0$. Then $V_{X_2}(-iA-lB)=V_{X_2}(-iA-(l+1)B)$. On the other hand, $d-l>c_r$, and hence $V_{X_3}(-(d-l)B)=0$. It follows that
\begin{equation*}
\begin{split}
\text{dim}\,K_{il}&=\text{dim}\,V_{X_1}(-(d-i)A)+\text{dim}\,V_{X_2}(-iA-lB)+\text{dim}\,V_{X_3}(-(d-l)B)\\
&=\text{dim}\,V_{X_2}(-iA-lB)\\
&=\text{dim}\,V_{X_2}=r+1.
\end{split}
\end{equation*}
We used above that $V_{X_2}(-iA)=V_{X_2}$ and $V_{X_2}(-lB)=V_{X_2}$.

Finally, suppose $l>b'_r$. Then 
\begin{center}
$V_{X_2}(-iA-lB)=V_{X_2}(-iA-(l+1)B)$ and $V_{X_3}(-(d-l)B)=V_{X_3}(-(d-l+1)B)$.
\end{center}
On the other hand, since $l>b'_r$, we have $d-l<c_0$, and hence $V_{X_3}(-(d-l)B)=V_{X_3}$. Then
\begin{equation*}
\begin{split}
\text{dim}\,K_{il}&=\text{dim}\,V_{X_1}(-(d-i)A)+\text{dim}\,V_{X_2}(-iA-lB)+\text{dim}\,V_{X_3}(-(d-l)B)\\
&=0+0+(r+1)=r+1,
\end{split}
\end{equation*}
We used above that $V_{X_2}(-lB)=0$, as $l>b'_r$. This finishes the proof of the stated inequality.

Now, we will prove that dim\,$K_{il}=r+1$ if $i\leq b_0$ or $l\leq b'_0$. We will only prove the stated equality in the case $i\leq b_0$, as the other case is analogous. Notice that, in Case 5, we saw dim\,$K_{il}=r+1$ if $i<b_0$. Thus, it remains to show the stated equality in the case $i=b_0$. Assume $i=b_0$.

Suppose first that $l=b'_k$ for some $k$. Notice that, in Case 1, for $j=0$, the equality holds in $\text{dim}\,V_{X_2}(-iA-lB)\geq r+1-(j+k)$, as $V_{X_2}(-iA-lB)=V_{X_2}(-lB)$. Thus dim\,$K_{il}=r+1$.

An analogous reasoning works for the case $b'_{k-1}<l<b'_k$. Now, suppose $l<b'_0$. In Case 5 we saw dim\,$K_{il}=r+1$ if $i<b_0$. Analogously, we can show that dim\,$K_{il}=r+1$ if $l<b'_0$. 

Finally, suppose that $l>b'_r$. Then 
\begin{center}
$V_{X_3}(-(d-l)B)=V_{X_3}(-(d-l+1)B)$, $V_{X_2}(-iA-lB)=0$ and $V_{X_3}(-(d-l)B)=V_{X_3}$. 
\end{center}
On the other hand, since $i=b_0$, we have $d-i=a_r$. It follows that
\begin{equation*}
\begin{split}
\text{dim}\,K_{il}&=\text{dim}\,V_{X_1}(-(d-i)A)+\text{dim}\,V_{X_2}(-iA-lB)+\text{dim}\,V_{X_3}(-(d-l)B)-1\\
&=(r+1-r)+0+(r+1)-1=r+1.
\end{split}
\end{equation*}
This finishes the proof of the stated equality.

Now, we will prove the statement 2 of the proposition. Keep the notation of multidegrees $\underline{d}$ and $\widetilde{\underline{d}}$ used in Section \ref{limit}. We will only prove the linking condition for $q=1$, as the proofs for $q=2,3$ are analogous. We will first prove that $\varphi_{\underline{d},\widetilde{\underline{d}}}(K_{il})\subseteq K_{i-1,l}$. (Recall that, for $q=1$, $\widetilde{\underline{d}}=(i-1,d-i-l+1,l)$.) Let $s\in K_{il}$. We have
\begin{equation}\label{link1}
(\alpha_{1 \widetilde{\underline{d}}}\circ \varphi_{\underline{d},\widetilde{\underline{d}}})(s)=(\alpha_{1 \underline{d}}\circ \varphi_{\widetilde{\underline{d}},\underline{d}}\circ \varphi_{\underline{d},\widetilde{\underline{d}}})(s)=(\alpha_{1 \underline{d}}\circ 0)(s)=0\in V_{X_1}(-(d-i+1)A).
\end{equation}
On the other hand, $s\in (\alpha_{2 \underline{d}})^{-1}(V_{X_2}(-iA-lB))\subseteq (\alpha_{2 \underline{d}})^{-1}(V_{X_2}(-(i-1)A-lB))$, as $s\in K_{il}$. Then
\begin{equation}\label{link2}
(\alpha_{2 \widetilde{\underline{d}}}\circ \varphi_{\underline{d},\widetilde{\underline{d}}})(s)=\alpha_{2 \underline{d}}(s)\in V_{X_2}(-(i-1)A-lB).
\end{equation}
Also, since $s\in K_{il}$, $s\in (\alpha_{3 \underline{d}})^{-1}(V_{X_3}(-(d-l)B))$, and hence
\begin{equation}\label{link3}
(\alpha_{3 \widetilde{\underline{d}}}\circ \varphi_{\underline{d},\widetilde{\underline{d}}})(s)=\alpha_{3 \underline{d}}(s)\in V_{X_3}(-(d-l)B).
\end{equation}
It follows from (\ref{link1}), (\ref{link2}) and (\ref{link3}) that $\varphi_{\underline{d},\widetilde{\underline{d}}}(s)\in K_{i-1,l}$. This proves that $\varphi_{\underline{d},\widetilde{\underline{d}}}(K_{il})\subseteq K_{i-1,l}$. 

Now, we will prove that $\varphi_{\widetilde{\underline{d}},\underline{d}}(K_{i-1,l})\subseteq K_{il}$. Let $s\in K_{i-1,l}$. Then 
\begin{center}
$s\in (\alpha_{1 \widetilde{\underline{d}}})^{-1}(V_{X_1}(-(d-i+1)A))\subseteq (\alpha_{1 \widetilde{\underline{d}}})^{-1}(V_{X_1}(-(d-i)A))$,
\end{center}
and hence
\begin{equation}\label{link4}
(\alpha_{1 \underline{d}}\circ \varphi_{\widetilde{\underline{d}},\underline{d}})(s)=\alpha_{1 \widetilde{\underline{d}}}(s)\in V_{X_1}(-(d-i)A).
\end{equation}
On the other hand, 
\begin{equation}\label{link5}
(\alpha_{2 \underline{d}}\circ \varphi_{\widetilde{\underline{d}},\underline{d}})(s)=(\alpha_{2 \widetilde{\underline{d}}}\circ \varphi_{\underline{d},\widetilde{\underline{d}}}\circ \varphi_{\widetilde{\underline{d}},\underline{d}})(s)=(\alpha_{2 \widetilde{\underline{d}}}\circ 0)(s)=0\in V_{X_2}(-iA-lB),
\end{equation}
and analogously
\begin{equation}\label{link6}
(\alpha_{3 \underline{d}}\circ \varphi_{\widetilde{\underline{d}},\underline{d}})(s)=(\alpha_{3 \widetilde{\underline{d}}}\circ \varphi_{\underline{d},\widetilde{\underline{d}}}\circ \varphi_{\widetilde{\underline{d}},\underline{d}})(s)=(\alpha_{3 \widetilde{\underline{d}}}\circ 0)(s)=0\in V_{X_3}(-(d-l)B).
\end{equation}
It follows from (\ref{link4}), (\ref{link5}) and (\ref{link6}) that $\varphi_{\widetilde{\underline{d}},\underline{d}}(s)\in K_{il}$. This proves that $\varphi_{\widetilde{\underline{d}},\underline{d}}(K_{i-1,l})\subseteq K_{il}$, which finishes the proof of the proposition.    \hfill $\Box$

\begin{prop}\label{ultrarefinado es exacto}
The following statements hold:

1. For any $i\geq 1$ and $l\geq 1$ such that $i+l\leq d$, 
\begin{center}
$\varphi_{\underline{d},\underline{d}'}(K_{il})=K_{i-1,l-1}^{X_{2}^{c},0}$,
\end{center}
where $\underline{d}:=(i,d-i-l,l)$ and $\underline{d}':=(i-1,d-i-l+2,l-1)$.

2. For any $i\geq 1$ and $l\geq 0$ such that $i+l\leq d$,
\begin{center}
$\varphi_{\underline{d},\underline{d}''}(K_{il})=K_{i-1,l}^{X_{1},0}$,
\end{center}
where $\underline{d}:=(i,d-i-l,l)$ and $\underline{d}'':=(i-1,d-i-l+1,l)$.

3. For any $i\geq 0$ and $l\geq 1$ such that $i+l\leq d$,
\begin{center}
$\varphi_{\underline{d},\widetilde{\underline{d}}}(K_{il})=K_{i,l-1}^{X_{3},0}$,
\end{center}
where $\underline{d}:=(i,d-i-l,l)$ and $\widetilde{\underline{d}}:=(i,d-i-l+1,l-1)$.
\end{prop}
{\em Proof.} We will first see how the statements 2 and 3 imply the statement 1. Let $i\geq 1$ and $l\geq 1$ such that $i+l\leq d$. Let $s'\in K_{i-1,l-1}^{X_{2}^{c},0}$. Then $s'\in K_{i-1,l-1}^{X_{3},0}$. But, by the statement 3 of the proposition, $\varphi_{\underline{d}'',{\underline{d}}'}(K_{i-1,l})=K_{i-1,l-1}^{X_{3},0}$, so $s'=\varphi_{\underline{d}'',{\underline{d}}'}(s'')$ for some $s''\in K_{i-1,l}$. As $\varphi_{\underline{d}'',{\underline{d}}'}(s'')=s'\in K_{i-1,l-1}^{X_{2}^{c},0}\subseteq K_{i-1,l-1}^{X_{1},0}$, it follows from Remark \ref{anulamiento} that $s''\in K_{i-1,l}^{X_{1},0}$. Then by the statement 2 of the proposition, $s''=\varphi_{\underline{d},{\underline{d}}''}(s)$ for some $s\in K_{il}$. Thus
\begin{center}
$s'=\varphi_{\underline{d}'',{\underline{d}}'}(s'')=\varphi_{\underline{d}'',{\underline{d}}'}\circ \varphi_{\underline{d},{\underline{d}}''}(s)=\varphi_{\underline{d},{\underline{d}}'}(s)\in \varphi_{\underline{d},\underline{d}'}(K_{il})$.
\end{center}
This proves that $K_{i-1,l-1}^{X_{2}^{c},0}\subseteq \varphi_{\underline{d},\underline{d}'}(K_{il})$. But, it follows from Proposition \ref{ultrarefinado es lls}, item 2, that $\varphi_{\underline{d},\underline{d}'}(K_{il})\subseteq K_{i-1,l-1}^{X_{2}^{c},0}$, so $\varphi_{\underline{d},\underline{d}'}(K_{il})=K_{i-1,l-1}^{X_{2}^{c},0}$.

It remains to show the statements 2 and 3. We will only prove the statement 2, as the statement 3 is analogous. 

By abuse of notation, we denote the restriction of $ev^{il}$ to the vector subspace 
\begin{center}
$V_{X_1}(-(d-i)A)\oplus V_{X_2}(-iA-lB)\oplus V_{X_3}(-(d-l)B)$
\end{center}
by $ev^{il}$ as well. It follows from the exact sequence defining $K_{il}$ that
\begin{equation*}
\text{dim}\,K_{il}=\text{dim}\,V_{X_1}(-(d-i)A)+\text{dim}\,V_{X_2}(-iA-lB)+\text{dim}\,V_{X_3}(-(d-l)B)-\text{dim}\,\text{Im}(ev^{il}).
\end{equation*}
On the other hand, the exact sequence defining $K_{il}$ induces the following exact sequence
\[0\rightarrow K_{il}^{X_{1}^{c},0}\rightarrow V_{X_1}(-(d-i)A)\oplus \{0\}\oplus \{0\}\rightarrow k\oplus k\]
Then $K_{il}^{X_{1}^{c},0}\cong  V_{X_1}(-(d-i+1)A)$, so dim\,$\varphi_{\underline{d},{\underline{d}}''}(K_{il})=$dim\,$K_{il}-$dim\,$V_{X_1}(-(d-i+1)A)$.
Thus
\begin{align}\label{casiexacto1}
\text{dim}\,\varphi_{\underline{d},{\underline{d}}''}(K_{il})=& \text{dim}\,V_{X_1}(-(d-i)A)+\text{dim}\,V_{X_2}(-iA-lB)+\text{dim}\,V_{X_3}(-(d-l)B) \nonumber \\
&-\text{dim}\,\text{Im}(ev^{il})-\text{dim}\,V_{X_1}(-(d-i+1)A). 
\end{align}
On the other hand, the exact sequence
\[0\rightarrow K_{i-1,l}^{X_1,0}\rightarrow \{0\}\oplus V_{X_2}(-(i-1)A-lB)\oplus V_{X_3}(-(d-l)B)\rightarrow k\oplus k\]
implies
\begin{equation}\label{casiexacto2}
\text{dim}\,K_{i-1,l}^{X_1,0}= \text{dim}\,V_{X_2}(-(i-1)A-lB)+\text{dim}\,V_{X_3}(-(d-l)B)-\text{dim}\,\text{Im}(ev_{1}^{i-1,l}),
\end{equation}
where $ev_{1}^{i-1,l}$ is the restriction of $ev^{i-1,l}$ to $\{0\}\oplus V_{X_2}(-(i-1)A-lB)\oplus V_{X_3}(-(d-l)B)$. But, it follows from Proposition \ref{ultrarefinado es lls}, item 2, that $\varphi_{\underline{d},{\underline{d}}''}(K_{il})\subseteq K_{i-1,l}^{X_1,0}$, so from (\ref{casiexacto1}) and (\ref{casiexacto2}), we have that $\varphi_{\underline{d},{\underline{d}}''}(K_{il})=K_{i-1,l}^{X_1,0}$ if and only if 
\begin{align}\label{casiexacto3}
&\text{dim}\,\text{Im}(ev_{1}^{i-1,l})-(\text{dim}\,V_{X_2}(-(i-1)A-lB)-\text{dim}\,V_{X_2}(-iA-lB)) \nonumber \\
&=\text{dim}\,\text{Im}(ev^{il})-(\text{dim}\,V_{X_1}(-(d-i)A)-\text{dim}\,V_{X_1}(-(d-i+1)A))
\end{align}
By checking cases $i=b_j$, $i\neq b_j$, $l=b'_k$ and $l\neq b'_k$, we see that both sides of (\ref{casiexacto3}) are equal to $\text{dim}\,V_{X_3}(-(d-l)B)-\text{dim}\,V_{X_3}(-(d-l+1)B)$. Thus (\ref{casiexacto3}) is true, and hence $\varphi_{\underline{d},{\underline{d}}''}(K_{il})=K_{i-1,l}^{X_1,0}$, proving the statement 2 of the proposition. This finishes the proof of the proposition.  \hfill $\Box$

\begin{prop}\label{exacto otro sentido}
The following statements hold:

1. For any $i\geq 1$ and $l\geq 0$ such that $i+l\leq d$, the following statements are equivalent:

(i) $\varphi_{\underline{d}'',\underline{d}}(K_{i-1,l})\neq K_{i,l}^{X_{1}^{c},0}$, where $\underline{d}:=(i,d-i-l,l)$ and $\underline{d}'':=(i-1,d-i-l+1,l)$.

(ii) $i-1$ is an order of vanishing of $V_{X_2}$ at $A$ and $i-1$ is not an order of vanishing of $V_{X_2}(-lB)$ at $A$.

2. For any $i\geq 0$ and $l\geq 1$ such that $i+l\leq d$, the following statements are equivalent:

(i) $\varphi_{\widetilde{\underline{d}},\underline{d}}(K_{i,l-1})\neq K_{i,l}^{X_{3}^{c},0}$, where $\underline{d}:=(i,d-i-l,l)$ and $\widetilde{\underline{d}}:=(i,d-i-l+1,l-1)$.

(ii) $l-1$ is an order of vanishing of $V_{X_2}$ at $B$ and $l-1$ is not an order of vanishing of $V_{X_2}(-iA)$ at $B$.
\end{prop}
{\em Proof.} We will only prove the statement 1. (The statement 2 is analogous.) Suppose first that $(i)$ holds. The exact sequence
\[0\rightarrow K_{il}^{X_{1}^{c},0}\rightarrow V_{X_1}(-(d-i)A)\oplus \{0\}\oplus \{0\}\rightarrow k\oplus k\]
implies that $\alpha_{1 \underline{d}}\big|_{K_{il}^{X_{1}^{c},0}}: K_{il}^{X_{1}^{c},0}\rightarrow V_{X_1}(-(d-i+1)A)$ is an isomorphism. On the other hand, it follows from Proposition \ref{ultrarefinado es lls}, item 2, that $\varphi_{\underline{d}'',{\underline{d}}}(K_{i-1,l})\subseteq K_{il}^{X_{1}^{c},0}$. Then, we have that $\varphi_{\underline{d}'',{\underline{d}}}(K_{i-1,l})=K_{il}^{X_{1}^{c},0}$ if and only if $\alpha_{1 \underline{d}}(\varphi_{\underline{d}'',{\underline{d}}}(K_{i-1,l}))=\alpha_{1 \underline{d}}(K_{il}^{X_{1}^{c},0})$, i.e., if and only if $\alpha_{1 \underline{d}''}(K_{i-1,l})=V_{X_1}(-(d-i+1)A)$. By hypothesis, $\varphi_{\underline{d}'',{\underline{d}}}(K_{i-1,l})\neq K_{il}^{X_{1}^{c},0}$, so $\alpha_{1 \underline{d}''}(K_{i-1,l})$ is a proper subspace of $V_{X_1}(-(d-i+1)A)$. The exact sequence defining $K_{i-1,l}$ induces the following exact sequence
\[0\rightarrow K_{i-1,l}\rightarrow \alpha_{1 \underline{d}''}(K_{i-1,l})\oplus V_{X_2}(-(i-1)A-lB)\oplus V_{X_3}(-(d-l)B)\rightarrow k\oplus k\]
By abuse of notation, we denote the restriction of $ev^{i-1,l}$ to the vector subspace 
\begin{center}
$V_{X_1}(-(d-i+1)A)\oplus V_{X_2}(-(i-1)A-lB)\oplus V_{X_3}(-(d-l)B)$
\end{center}
by $ev^{i-1,l}$ as well, and let $\overline{ev}^{i-1,l}$ be the restriction of $ev^{i-1,l}$ to the vector subspace $\alpha_{1 \underline{d}''}(K_{i-1,l})\oplus V_{X_2}(-(i-1)A-lB)\oplus V_{X_3}(-(d-l)B)$. We have 
\begin{equation*}
\begin{split}
\text{dim}\,K_{i-1,l}= & \text{dim}\,V_{X_1}(-(d-i+1)A)+\text{dim}\,V_{X_2}(-(i-1)A-lB)\\
&+\text{dim}\,V_{X_3}(-(d-l)B)-\text{dim}\,\text{Im}(ev^{i-1,l})
\end{split}
\end{equation*}
and also
\begin{equation*}
\begin{split}
\text{dim}\,K_{i-1,l}= & \text{dim}\,\alpha_{1 \underline{d}''}(K_{i-1,l})+\text{dim}\,V_{X_2}(-(i-1)A-lB)+\text{dim}\,V_{X_3}(-(d-l)B)\\
&-\text{dim}\,\text{Im}(\overline{ev}^{i-1,l}).
\end{split}
\end{equation*}
Therefore
\begin{equation}\label{imagen1}
\begin{split}
\text{dim}\,V_{X_1}(-(d-i+1)A)-\text{dim}\,\alpha_{1 \underline{d}''}(K_{i-1,l})=\text{dim}\,\text{Im}(ev^{i-1,l})-\text{dim}\,\text{Im}(\overline{ev}^{i-1,l}),
\end{split}
\end{equation}
and since $\text{dim}\,\text{Im}(ev^{i-1,l})-\text{dim}\,\text{Im}(\overline{ev}^{i-1,l})\leq \text{dim}\,\text{Im}(ev^{i-1,l})\leq 2$, it follows that
\begin{center}
$\text{dim}\,V_{X_1}(-(d-i+1)A)-2\leq \text{dim}\,\alpha_{1 \underline{d}''}(K_{i-1,l})\leq \text{dim}\,V_{X_1}(-(d-i+1)A)-1$, 
\end{center}
as $\alpha_{1 \underline{d}''}(K_{i-1,l})$ is a proper subspace of $V_{X_1}(-(d-i+1)A)$. Thus, there are two cases to consider.\\
{\em Case 1:} If $\text{dim}\,\alpha_{1 \underline{d}''}(K_{i-1,l})=\text{dim}\,V_{X_1}(-(d-i+1)A)-1$.

It follows from (\ref{imagen1}) that $\text{dim}\,\text{Im}(\overline{ev}^{i-1,l})=\text{dim}\,\text{Im}(ev^{i-1,l})-1$. We will first prove that $i-1=b_j$ for some $j$. Suppose by contradiction that $i-1$ is not an order of vanishing of $V_{X_2}$ at $A$. Then $\text{dim}\,\text{Im}(ev^{i-1,l})\leq 1$, and hence $\text{dim}\,\text{Im}(\overline{ev}^{i-1,l})\leq 0$. So $\text{dim}\,\text{Im}(\overline{ev}^{i-1,l})=0$ and $\text{dim}\,\text{Im}(ev^{i-1,l})=1$. Now, since $i-1$ is not an order of vanishing of $V_{X_2}$ at $A$, $\text{dim}\,\text{Im}(ev^{i-1,l})=1$ implies that $l=b'_k$ for some $k$, which implies that $\text{Im}(\overline{ev}^{i-1,l})\supseteq \{0\}\oplus k$, and hence $\text{dim}\,\text{Im}(\overline{ev}^{i-1,l})\geq 1$, a contradiction. Thus $i-1=b_j$ for some $j$. 

Now, we will prove that $i-1$ is not an order of vanishing of $V_{X_2}(-lB)$ at $A$. Suppose first that $l=b'_k$ for some $k$. Since $i-1=b_j$ and $l=b'_k$, we have $\text{dim}\,\text{Im}(ev^{i-1,l})=2$. Then $\text{dim}\,\text{Im}(\overline{ev}^{i-1,l})=1$. Since $l=b'_k$, $\text{Im}(\overline{ev}^{i-1,l})\supseteq \{0\}\oplus k$, and it follows from dimension considerations that $\text{Im}(\overline{ev}^{i-1,l})=\{0\}\oplus k$. This implies that all sections of $\alpha_{1 \underline{d}''}(K_{i-1,l})\subseteq V_{X_1}(-(d-i+1)A)$ and all sections of $V_{X_2}(-(i-1)A-lB)$ vanish at $A$, i.e., 
\begin{center}
$\alpha_{1 \underline{d}''}(K_{i-1,l})\subseteq V_{X_1}(-(d-i+2)A)$ and 
$V_{X_2}(-(i-1)A-lB)=V_{X_2}(-iA-lB)$. 
\end{center}
Thus $i-1$ is not an order of vanishing of $V_{X_2}(-lB)$ at $A$. In addition, since $\alpha_{1 \underline{d}''}(K_{i-1,l})$ is contained in $V_{X_1}(-(d-i+2)A)$, it follows from dimension considerations that $\alpha_{1 \underline{d}''}(K_{i-1,l})$ is equal to $V_{X_1}(-(d-i+2)A)$.

Now, assume $l$ is not an order of vanishing of $V_{X_2}$ at $B$. Then $\text{dim}\,\text{Im}(ev^{i-1,l})=1$, as $i-1=b_j$, and hence $\text{Im}(\overline{ev}^{i-1,l})=\{0\}\oplus \{0\}$. This implies that all sections of $\alpha_{1 \underline{d}''}(K_{i-1,l})\subseteq V_{X_1}(-(d-i+1)A)$ and all sections of $V_{X_2}(-(i-1)A-lB)$ vanish at $A$. It follows that $\alpha_{1 \underline{d}''}(K_{i-1,l})=V_{X_1}(-(d-i+2)A)$ and $i-1$ is not an order of vanishing of $V_{X_2}(-lB)$ at $A$.\\
{\em Case 2:} If $\text{dim}\,\alpha_{1 \underline{d}''}(K_{i-1,l})=\text{dim}\,V_{X_1}(-(d-i+1)A)-2$.

It follows from (\ref{imagen1}) that $\text{dim}\,\text{Im}(\overline{ev}^{i-1,l})=\text{dim}\,\text{Im}(ev^{i-1,l})-2$. As $\text{dim}\,\text{Im}(ev^{i-1,l})\leq 2$, we have $\text{Im}(\overline{ev}^{i-1,l})=\{0\}\oplus \{0\}$ and $\text{dim}\,\text{Im}(ev^{i-1,l})=2$. Since $\text{dim}\,\text{Im}(ev^{i-1,l})=2$, it follows that $i-1=b_j$ and $l=b'_k$ for some $j,k$. Now, as $l=b'_k$, we get $\text{dim}\,\text{Im}(\overline{ev}^{i-1,l})\geq 1$, a contradiction. Thus, the only case can happen is Case 1, and hence $(ii)$ holds.  

Suppose now that $(ii)$ holds. Define $K'\subseteq H^{0}(\mathcal{L}_{\underline{d}''})$ by the exact sequence
\[0\rightarrow K'\rightarrow V_{X_1}(-(d-i+2)A)\oplus V_{X_2}(-(i-1)A-lB)\oplus V_{X_3}(-(d-l)B)\rightarrow k\oplus k.\]
By abuse of notation, we denote the restriction of $ev^{i-1,l}$ to the vector subspace 
\begin{center}
$V_{X_1}(-(d-i+1)A)\oplus V_{X_2}(-(i-1)A-lB)\oplus V_{X_3}(-(d-l)B)$
\end{center}
by $ev^{i-1,l}$ as well, and let $\widetilde{ev}^{i-1,l}$ be the restriction of $ev^{i-1,l}$ to the vector subspace $V_{X_1}(-(d-i+2)A)\oplus V_{X_2}(-(i-1)A-lB)\oplus V_{X_3}(-(d-l)B)$. We have that $K'\subseteq K_{i-1,l}$, as $V_{X_1}(-(d-i+2)A)\subseteq V_{X_1}(-(d-i+1)A)$. On the other hand, it follows from the definition of $K'$ that $\alpha_{1 \underline{d}''}(K')\subseteq V_{X_1}(-(d-i+2)A)$, and hence $\alpha_{1 \underline{d}''}(K')\neq V_{X_1}(-(d-i+1)A)$. Now, recall that $\varphi_{\underline{d}'',{\underline{d}}}(K_{i-1,l})=K_{il}^{X_{1}^{c},0}$ if and only if $\alpha_{1 \underline{d}''}(K_{i-1,l})=V_{X_1}(-(d-i+1)A)$. Therefore, to prove $(i)$, we need only show that $\alpha_{1 \underline{d}''}(K_{i-1,l})\neq V_{X_1}(-(d-i+1)A)$. For this, it suffices to show that $K_{i-1,l}=K'$.

Since $K'\subseteq K_{i-1,l}$, we need only prove that dim\,$K'=$dim\,$K_{i-1,l}$. We have 
\begin{equation*}
\begin{split}
\text{dim}\,K_{i-1,l}= & \text{dim}\,V_{X_1}(-(d-i+1)A)+\text{dim}\,V_{X_2}(-(i-1)A-lB)\\
&+\text{dim}\,V_{X_3}(-(d-l)B)-\text{dim}\,\text{Im}(ev^{i-1,l})
\end{split}
\end{equation*}
and also
\begin{equation*}
\begin{split}
\text{dim}\,K'= & \text{dim}\,V_{X_1}(-(d-i+2)A)+\text{dim}\,V_{X_2}(-(i-1)A-lB)\\
&+\text{dim}\,V_{X_3}(-(d-l)B)-\text{dim}\,\text{Im}(\widetilde{ev}^{i-1,l}).
\end{split}
\end{equation*}
Therefore dim\,$K'=$dim\,$K_{i-1,l}$ if and only if
\begin{equation*}
\text{dim}\,V_{X_1}(-(d-i+1)A)-\text{dim}\,V_{X_1}(-(d-i+2)A)=\text{dim}\,\text{Im}(ev^{i-1,l})-\text{dim}\,\text{Im}(\widetilde{ev}^{i-1,l}),
\end{equation*}
i.e., if and only if $\text{dim}\,\text{Im}(\widetilde{ev}^{i-1,l})=\text{dim}\,\text{Im}(ev^{i-1,l})-1$, as $i-1=b_j$ for some $j$. There are two cases to consider.\\
{\em Case I:} If $l=b'_k$ for some $k$.

Since $i-1=b_j$ and $l=b'_k$, $\text{dim}\,\text{Im}(ev^{i-1,l})=2$. On the other hand, since $i-1$ is not an order of vanishing of $V_{X_2}(-lB)$ at $A$, we have $V_{X_2}(-(i-1)A-lB)=V_{X_2}(-iA-lB)$, i.e., all sections of $V_{X_2}(-(i-1)A-lB)$ vanish at $A$. Then $\text{Im}(\widetilde{ev}^{i-1,l})\subseteq \{0\}\oplus k$. But, since $l=b'_k$, $\text{Im}(\widetilde{ev}^{i-1,l})\supseteq \{0\}\oplus k$, and hence $\text{Im}(\widetilde{ev}^{i-1,l})=\{0\}\oplus k$. Therefore $\text{dim}\,\text{Im}(\widetilde{ev}^{i-1,l})=1$, and thus $\text{dim}\,\text{Im}(\widetilde{ev}^{i-1,l})=\text{dim}\,\text{Im}(ev^{i-1,l})-1$.\\
{\em Case II:} If $l$ is not an order of vanishing of $V_{X_2}$ at $B$.

Since $i-1=b_j$, $\text{dim}\,\text{Im}(ev^{i-1,l})=1$. As in Case I, $\text{Im}(\widetilde{ev}^{i-1,l})\subseteq \{0\}\oplus k$. But, since $l$ is not an order of vanishing of $V_{X_2}$ at $B$, $\text{Im}(\widetilde{ev}^{i-1,l})\subseteq k\oplus \{0\}$, and hence we have $\text{Im}(\widetilde{ev}^{i-1,l})=\{0\}\oplus \{0\}$. So $\text{dim}\,\text{Im}(\widetilde{ev}^{i-1,l})=0=\text{dim}\,\text{Im}(ev^{i-1,l})-1$. This finishes the proof of the proposition.   \hfill $\Box$

\begin{rem}\label{desdoble}
{\em Let $V_1, V_2$ and $V_3$ vector subspaces of a $N$-dimensional vector space $V$. We will say that $V_1$ {\em distributes over} $V_2$ and $V_3$ if $V_1\cap (V_2+V_3)=V_1\cap V_2+V_1\cap V_3$. Note that this notion is symmetric on $V_1, V_2$ and $V_3$.

Indeed, since $V_1\cap V_2+V_1\cap V_3 \subseteq V_1\cap(V_2+V_3)$, we have that $V_1\cap (V_2+V_3)=V_1\cap V_2+V_1\cap V_3$ is equivalent to dim\,$V_1\cap (V_2+V_3)=$dim\,$(V_1\cap V_2+V_1\cap V_3)$. We have
\begin{center}
$\text{dim}\,V_1\cap (V_2+V_3)=\text{dim}\,(V_1\cap V_2+V_1\cap V_3)$ if and only if
\end{center}
$\text{dim}\,V_1+\text{dim}\,(V_2+V_3)-\text{dim}\,(V_1+V_2+V_3)=\text{dim}\,V_1\cap V_2+\text{dim}\,V_1\cap V_3-\text{dim}\,V_1\cap V_2\cap V_3$.\\
This is equivalent to
\begin{equation*}
\begin{split}
&\text{dim}\,V_1+\text{dim}\,V_2+\text{dim}\,V_3-\text{dim}\,V_2\cap V_3-\text{dim}\,(V_1+V_2+V_3)\\
&=\text{dim}\,V_1\cap V_2+\text{dim}\,V_1\cap V_3-\text{dim}\,V_1\cap V_2\cap V_3, \text{i.e.},
\end{split}
\end{equation*}
\begin{align}\label{desdoble1}
\text{dim}\,(V_1+V_2+V_3)=&\text{dim}\,V_1+\text{dim}\,V_2+\text{dim}\,V_3-\text{dim}\,V_1\cap V_2-\text{dim}\,V_1\cap V_3-\text{dim}\,V_2\cap V_3 \nonumber \\
&+\text{dim}\,V_1\cap V_2\cap V_3.
\end{align}
Now, just notice that (\ref{desdoble1}) is symmetric on $V_1, V_2$ and $V_3$. Thus, the following statements are equivalent:\\
1. $V_1$ distributes over $V_2$ and $V_3$.\\
2. $V_2$ distributes over $V_1$ and $V_3$.\\
3. $V_3$ distributes over $V_1$ and $V_2$.
}
\end{rem}

\begin{prop}\label{desdoblekernel}
For any $i\geq 0$ and $l\geq 0$ such that $i+l\leq d$:

$K_{il}^{X_{q_1},0}\cap (K_{il}^{X_{q_2},0}+K_{il}^{X_{q_3},0})=K_{il}^{X_{q_3}^{c},0}+K_{il}^{X_{q_2}^{c},0}$ for distinct $q_1,q_2,q_3$.


\end{prop}
{\em Proof.} By Remark \ref{desdoble}, it is enough to prove the case $q_{m}=m$ for $m=1,2,3$. Via the injective map in the exact sequence defining $K_{il}$, we can see $K_{il}$ as a subspace of $V_{X_1}(-(d-i)A)\oplus V_{X_2}(-iA-lB)\oplus V_{X_3}(-(d-l)B)$. It follows from the exact sequence defining $K_{il}$ that
\begin{center}
$K_{il}^{X_2,0}=V_{X_1}(-(d-i+1)A)\oplus \{0\}\oplus V_{X_3}(-(d-l+1)B)$,\\
$K_{il}^{X_{3}^{c},0}=\{0\}\oplus \{0\}\oplus V_{X_3}(-(d-l+1)B)$,\\
$K_{il}^{X_{2}^{c},0}=\{0\}\oplus V_{X_2}(-(i+1)A-(l+1)B)\oplus \{0\}$\\
$\text{and}\,\,\,K_{il}^{X_{1}^{c},0}=V_{X_1}(-(d-i+1)A)\oplus \{0\}\oplus \{0\}$.\\
\end{center}
Also, by checking cases $i=b_j$, $i\neq b_j$, $l=b'_k$ and $l\neq b'_k$, we get
\begin{center}
$\text{dim}\,K_{il}^{X_1,0}=\text{dim}\,V_{X_2}(-(i+1)A-lB)+\text{dim}\,V_{X_3}(-(d-l+1)B)$,\\
$\text{dim}\,K_{il}^{X_3,0}=\text{dim}\,V_{X_1}(-(d-i+1)A)+\text{dim}\,V_{X_2}(-iA-(l+1)B)$\\
$\text{and}\,\,\,\text{dim}\,K_{il}=\text{dim}\,V_{X_1}(-(d-i+1)A)+\text{dim}\,V_{X_2}(-iA-lB)+\text{dim}\,V_{X_3}(-(d-l+1)B).$
\end{center}
Then
\begin{align} 
\text{dim}(K_{il}^{X_2,0}+K_{il}^{X_3,0})=& \text{dim}\,V_{X_1}(-(d-i+1)A)+\text{dim}\,V_{X_2}(-iA-(l+1)B) \nonumber \\
& +\text{dim}\,V_{X_3}(-(d-l+1)B). \nonumber
\end{align}
In particular, $\text{dim}(K_{il}^{X_2,0}+K_{il}^{X_3,0})\geq \text{dim}\,K_{il}-1$, and equality holds if and only if $l$ is an order of vanishing of $V_{X_2}(-iA)$ at $B$. On the other hand, we have
\begin{center}
$\text{dim}(K_{il}^{X_{3}^{c},0}+K_{il}^{X_{2}^{c},0})=\text{dim}\,V_{X_2}(-(i+1)A-(l+1)B)+\text{dim}\,V_{X_3}(-(d-l+1)B)$, 
\end{center}
so, by dimension considerations, we have that $K_{il}^{X_{3}^{c},0}+K_{il}^{X_{2}^{c},0}=K_{il}^{X_1,0}$ if and only if $V_{X_2}(-(i+1)A-lB)=V_{X_2}(-(i+1)A-(l+1)B)$. In this case, the statement of the proposition holds. 

Suppose now that $V_{X_2}(-(i+1)A-lB)\neq V_{X_2}(-(i+1)A-(l+1)B)$. Then, we have $V_{X_2}(-iA-lB)\neq V_{X_2}(-iA-(l+1)B)$, and hence $\text{dim}(K_{il}^{X_2,0}+K_{il}^{X_3,0})=\text{dim}\,K_{il}-1$. Notice that
\begin{align}
\text{dim}(K_{il}^{X_{3}^{c},0}+K_{il}^{X_{2}^{c},0})&=\text{dim}\,V_{X_2}(-(i+1)A-(l+1)B)+\text{dim}\,V_{X_3}(-(d-l+1)B) \nonumber \\
&= \text{dim}\,V_{X_2}(-(i+1)A-lB)-1+\text{dim}\,V_{X_3}(-(d-l+1)B) \nonumber \\
&=\text{dim}\,K_{il}^{X_1,0}-1. \nonumber
\end{align}
Now, since $K_{il}^{X_{3}^{c},0}+K_{il}^{X_{2}^{c},0}\subseteq K_{il}^{X_1,0}\cap (K_{il}^{X_2,0}+K_{il}^{X_3,0})$ and $\text{dim}(K_{il}^{X_2,0}+K_{il}^{X_3,0})=\text{dim}\,K_{il}-1$, it suffices to show that $K_{il}^{X_1,0}$ is not contained in the space $K_{il}^{X_2,0}+K_{il}^{X_3,0}$. Suppose by contradiction that $K_{il}^{X_1,0}\subseteq K_{il}^{X_2,0}+K_{il}^{X_3,0}$. Notice that 
\begin{center}
$K_{il}^{X_3,0}\subseteq V_{X_1}(-(d-i)A)\oplus V_{X_2}(-iA-(l+1)B)\oplus \{0\}$. 
\end{center}
It follows that
\begin{center}
$K_{il}^{X_2,0}+K_{il}^{X_3,0}\subseteq V_{X_1}(-(d-i)A)\oplus V_{X_2}(-iA-(l+1)B)\oplus V_{X_3}(-(d-l+1)B)$, 
\end{center}
and hence $K_{il}^{X_1,0}\subseteq V_{X_1}(-(d-i)A)\oplus V_{X_2}(-iA-(l+1)B)\oplus V_{X_3}(-(d-l+1)B)$. So, since $K_{il}^{X_1,0}\subseteq \{0\}\oplus V_{X_2}(-(i+1)A-lB)\oplus V_{X_3}(-(d-l)B)$, we get
\begin{center}
$K_{il}^{X_1,0}\subseteq \{0\}\oplus V_{X_2}(-(i+1)A-(l+1)B)\oplus V_{X_3}(-(d-l+1)B)$.
\end{center}
Then 
\begin{align}
\text{dim}\,K_{il}^{X_1,0}&\leq \text{dim}\,V_{X_2}(-(i+1)A-(l+1)B)+\text{dim}\,V_{X_3}(-(d-l+1)B) \nonumber \\
&=\text{dim}\,V_{X_2}(-(i+1)A-lB)-1+\text{dim}\,V_{X_3}(-(d-l+1)B) \nonumber \\
&=\text{dim}\,K_{il}^{X_1,0}-1, \nonumber
\end{align}
a contradiction. So the statement of the proposition is shown.   \hfill $\Box$

\begin{prop}\label{producto}
For any $i\geq 0$ and $l\geq 0$ such that $i+l\leq d$, the following statements hold:

1. {\em dim}$(K_{il}^{X_{q_{1}},0}+K_{il}^{X_{q_{2}},0})\geq$ {\em dim\,}$K_{il}-1$ for any $q_{1}\neq q_{2}$.



2. {\em dim\,}$K_{il}^{X_1,0}+${\em dim\,}$K_{il}^{X_2,0}+${\em dim\,}$K_{il}^{X_3,0}\geq 2(${\em dim\,}$K_{il}-1)$.
\end{prop}
{\em Proof.} By the proof of Proposition \ref{desdoblekernel}, the statement 1 holds for $q_{1}=2$ and $q_{2}=3$. The proofs of the other cases are analogous. As for the statement 2, putting together the following equalities
\begin{center}
$K_{il}^{X_2,0}=V_{X_1}(-(d-i+1)A)\oplus \{0\}\oplus V_{X_3}(-(d-l+1)B)$,\\
$\text{dim}\,K_{il}^{X_1,0}=\text{dim}\,V_{X_2}(-(i+1)A-lB)+\text{dim}\,V_{X_3}(-(d-l+1)B)$,\\
$\text{dim}\,K_{il}^{X_3,0}=\text{dim}\,V_{X_1}(-(d-i+1)A)+\text{dim}\,V_{X_2}(-iA-(l+1)B)$\\
$\text{and}\,\,\,\text{dim}\,K_{il}=\text{dim}\,V_{X_1}(-(d-i+1)A)+\text{dim}\,V_{X_2}(-iA-lB)+\text{dim}\,V_{X_3}(-(d-l+1)B),$
\end{center}
and the inequalities
\begin{align}
\text{dim}\,V_{X_2}(-(i+1)A-lB)\geq \text{dim}\,V_{X_2}(-iA-lB)-1 \nonumber\\
\text{and}\,\,\,\text{dim}\,V_{X_2}(-iA-(l+1)B)\geq \text{dim}\,V_{X_2}(-iA-lB)-1, \nonumber
\end{align}
we get dim\,$K_{il}^{X_1,0}+$dim\,$K_{il}^{X_2,0}+$dim\,$K_{il}^{X_3,0}\geq 2($dim\,$K_{il}-1)$, and equality holds if and only if the equality holds in the two inequalities above.   \hfill $\Box$

\begin{prop}\label{novacio} 
For any $i\geq 0$ and $l\geq 0$ such that $i+l\leq d$:

{\em dim\,}$K_{il}^{X_1,0}+${\em dim\,}$K_{il}^{X_2,0}+${\em dim\,}$K_{il}^{X_3,0}=2(${\em dim\,}$K_{il}-1)$ if and only if $i$ is an order of vanishing of $V_{X_2}(-lB)$ at $A$ and $l$ is an order of vanishing of $V_{X_2}(-iA)$ at $B$. In this case, we have that $K_{il}^{X_1,0}, K_{il}^{X_2,0}$ and $K_{il}^{X_3,0}$ are proper subspaces of $K_{il}$.
\end{prop}
{\em Proof.} We have seen in the proof of Proposition \ref{producto} that 
\begin{center}
$\text{dim}\,K_{il}^{X_1,0}+\text{dim}\,K_{il}^{X_2,0}+\text{dim}\,K_{il}^{X_3,0}=2(\text{dim}\,K_{il}-1)$ 
\end{center}
if and only if 
\begin{align}
\text{dim}\,V_{X_2}(-(i+1)A-lB)=\text{dim}\,V_{X_2}(-iA-lB)-1 \nonumber\\
\text{and}\,\,\,\text{dim}\,V_{X_2}(-iA-(l+1)B)=\text{dim}\,V_{X_2}(-iA-lB)-1, \nonumber
\end{align}
i.e., if and only if $i$ is an order of vanishing of $V_{X_2}(-lB)$ at $A$ and $l$ is an order of vanishing of $V_{X_2}(-iA)$ at $B$.

Now, suppose that $i$ is an order of vanishing of $V_{X_2}(-lB)$ at $A$ and $l$ is an order of vanishing of $V_{X_2}(-iA)$ at $B$. By the proof of Proposition \ref{desdoblekernel}, $K_{il}^{X_2,0}+K_{il}^{X_3,0}$ has dimension $\text{dim}\,K_{il}-1$ if and only if $l$ is an order of vanishing of $V_{X_2}(-iA)$ at $B$. So, by the hypothesis on $l$, we have $\text{dim}(K_{il}^{X_2,0}+K_{il}^{X_3,0})=\text{dim}\,K_{il}-1$, and hence $K_{il}^{X_2,0}$ and $K_{il}^{X_3,0}$ are proper subspaces of $K_{il}$. Analogously, since $i$ is an order of vanishing of $V_{X_2}(-lB)$ at $A$, we have $\text{dim}(K_{il}^{X_2,0}+K_{il}^{X_1,0})=\text{dim}\,K_{il}-1$, and hence $K_{il}^{X_1,0}$ is a proper subspace of $K_{il}$.   \hfill $\Box$

\section{Constructing exact extensions}\label{extension}

We will describe a method for the construction of exact extensions. Furthermore, this method allows us to construct any exact extension. The main result of this section is the following proposition, which is the fundamental statement for our method.

\begin{prop}\label{exact ext}
For any $i\geq 0$ and $l\geq 0$ such that $i+l\leq d$, let $\underline{d}:=(i,d-i-l,l)$, and let $\underline{d}'':=(i-1,d-i-l+1,l)$ if $i>0$. Then, the following statements hold:

1. Let $i,l$ be positive integers such that $i+l=d$.\\
Let $\underline{d}':=(i-1,d-i-l+2,l-1)$. Let $V_{\underline{d}'}$ and $V_{\underline{d}''}$ be $r+1$-dimensional subspaces of $K_{i-1,l-1}$ and $K_{i-1,l}$, respectively, such that 
\begin{center}
$\varphi_{\underline{d}',\underline{d}''}(V_{\underline{d}'})=V_{\underline{d}''}^{X_{3}^{c},0}$ and $\varphi_{\underline{d}'',\underline{d}'}(V_{\underline{d}''})=V_{\underline{d}'}^{X_{3},0}$.
\end{center}
Set $\beta:=${\em dim\,}$V_{\underline{d}''}^{X_{1},0}-${\em dim}$(V_{\underline{d}''}^{X_{2}^{c},0}\oplus V_{\underline{d}''}^{X_{3}^{c},0})$. Then, for any linearly independent elements $u_1,\ldots,u_{\beta}\in V_{\underline{d}''}^{X_{1},0}$ such that $V_{\underline{d}''}^{X_{1},0}=(V_{\underline{d}''}^{X_{2}^{c},0}\oplus V_{\underline{d}''}^{X_{3}^{c},0})\oplus \left\langle u_1,\ldots,u_{\beta}\right\rangle$, and for any elements $v_1,\ldots,v_{\beta}\in K_{il}$ such that $\varphi_{\underline{d},\underline{d}''}(v_1)=u_1,\ldots,\varphi_{\underline{d},\underline{d}''}(v_{\beta})=u_{\beta}$, the subspace 
\begin{center}
$V_{\underline{d}}:=\varphi_{\underline{d}',\underline{d}}(V_{\underline{d}'})+\left\langle v_1,\ldots,v_{\beta}\right\rangle \subseteq K_{il}$
\end{center}
is $r+1$-dimensional, and
\begin{center}
$\varphi_{\underline{d}',\underline{d}}(V_{\underline{d}'})=V_{\underline{d}}^{X_{2},0}$, $\varphi_{\underline{d},\underline{d}'}(V_{\underline{d}})=V_{\underline{d}'}^{X_{2}^{c},0}$ and\\
$\varphi_{\underline{d}'',\underline{d}}(V_{\underline{d}''})=V_{\underline{d}}^{X_{1}^{c},0}$, $\varphi_{\underline{d},\underline{d}''}(V_{\underline{d}})=V_{\underline{d}''}^{X_{1},0}$.
\end{center}

2. Let $i,l$ be positive integers such that $i+l\leq d-1$.\\ 
Let $\underline{d}':=(i-1,d-i-l+2,l-1)$ and $\underline{d}''':=(i,d-i-l-1,l+1)$. Let $V_{\underline{d}'}$, $V_{\underline{d}''}$ and $V_{\underline{d}'''}$ be $r+1$-dimensional subspaces of $K_{i-1,l-1}$, $K_{i-1,l}$ and $K_{i,l+1}$, respectively, such that 
\begin{center}
$\varphi_{\underline{d}',\underline{d}''}(V_{\underline{d}'})=V_{\underline{d}''}^{X_{3}^{c},0}$, $\varphi_{\underline{d}'',\underline{d}'}(V_{\underline{d}''})=V_{\underline{d}'}^{X_{3},0}$, $\varphi_{\underline{d}'',\underline{d}'''}(V_{\underline{d}''})=V_{\underline{d}'''}^{X_{2},0}$ and $\varphi_{\underline{d}''',\underline{d}''}(V_{\underline{d}'''})=V_{\underline{d}''}^{X_{2}^{c},0}$.
\end{center}
Set $\beta:=${\em dim\,}$V_{\underline{d}''}^{X_{1},0}-${\em dim}$(V_{\underline{d}''}^{X_{2}^{c},0}\oplus V_{\underline{d}''}^{X_{3}^{c},0})$. Then, for any linearly independent elements $u_1,\ldots,u_{\beta}\in V_{\underline{d}''}^{X_{1},0}$ such that $V_{\underline{d}''}^{X_{1},0}=(V_{\underline{d}''}^{X_{2}^{c},0}\oplus V_{\underline{d}''}^{X_{3}^{c},0})\oplus \left\langle u_1,\ldots,u_{\beta}\right\rangle$, and for any elements $v_1,\ldots,v_{\beta}\in K_{il}$ such that $\varphi_{\underline{d},\underline{d}''}(v_1)=u_1,\ldots,\varphi_{\underline{d},\underline{d}''}(v_{\beta})=u_{\beta}$, the subspace 
\begin{center}
$V_{\underline{d}}:=(\varphi_{\underline{d}',\underline{d}}(V_{\underline{d}'})+\varphi_{\underline{d}''',\underline{d}}(V_{\underline{d}'''}))+\left\langle v_1,\ldots,v_{\beta}\right\rangle \subseteq K_{il}$
\end{center}
is $r+1$-dimensional, and 
\begin{center}
$\varphi_{\underline{d}',\underline{d}}(V_{\underline{d}'})=V_{\underline{d}}^{X_{2},0}$, $\varphi_{\underline{d},\underline{d}'}(V_{\underline{d}})=V_{\underline{d}'}^{X_{2}^{c},0}$,\\
$\varphi_{\underline{d}'',\underline{d}}(V_{\underline{d}''})=V_{\underline{d}}^{X_{1}^{c},0}$, $\varphi_{\underline{d},\underline{d}''}(V_{\underline{d}})=V_{\underline{d}''}^{X_{1},0}$ and\\
$\varphi_{\underline{d}''',\underline{d}}(V_{\underline{d}'''})=V_{\underline{d}}^{X_{3},0}$, $\varphi_{\underline{d},\underline{d}'''}(V_{\underline{d}})=V_{\underline{d}'''}^{X_{3}^{c},0}$.
\end{center}

3. Let $0<i<d$ and $l=0$.\\
Let $\underline{d}''':=(i,d-i-l-1,l+1)$. Let $V_{\underline{d}''}$ and $V_{\underline{d}'''}$ be $r+1$-dimensional subspaces of $K_{i-1,l}$ and $K_{i,l+1}$, respectively, such that 
\begin{center}
$\varphi_{\underline{d}'',\underline{d}'''}(V_{\underline{d}''})=V_{\underline{d}'''}^{X_{2},0}$ and $\varphi_{\underline{d}''',\underline{d}''}(V_{\underline{d}'''})=V_{\underline{d}''}^{X_{2}^{c},0}$.
\end{center}
Set $\beta:=${\em dim\,}$V_{\underline{d}''}^{X_{1},0}-${\em dim}$(V_{\underline{d}''}^{X_{2}^{c},0}\oplus V_{\underline{d}''}^{X_{3}^{c},0})$. Then, for any linearly independent elements $u_1,\ldots,u_{\beta}\in V_{\underline{d}''}^{X_{1},0}$ such that $V_{\underline{d}''}^{X_{1},0}=(V_{\underline{d}''}^{X_{2}^{c},0}\oplus V_{\underline{d}''}^{X_{3}^{c},0})\oplus \left\langle u_1,\ldots,u_{\beta}\right\rangle$, and for any elements $v_1,\ldots,v_{\beta}\in K_{il}$ such that $\varphi_{\underline{d},\underline{d}''}(v_1)=u_1,\ldots,\varphi_{\underline{d},\underline{d}''}(v_{\beta})=u_{\beta}$, the subspace 
\begin{center}
$V_{\underline{d}}:=\varphi_{\underline{d}''',\underline{d}}(V_{\underline{d}'''})+\left\langle v_1,\ldots,v_{\beta}\right\rangle \subseteq K_{il}$
\end{center}
is $r+1$-dimensional, and
\begin{center}
$\varphi_{\underline{d}''',\underline{d}}(V_{\underline{d}'''})=V_{\underline{d}}^{X_{3},0}$, $\varphi_{\underline{d},\underline{d}'''}(V_{\underline{d}})=V_{\underline{d}'''}^{X_{3}^{c},0}$ and\\
$\varphi_{\underline{d}'',\underline{d}}(V_{\underline{d}''})=V_{\underline{d}}^{X_{1}^{c},0}$, $\varphi_{\underline{d},\underline{d}''}(V_{\underline{d}})=V_{\underline{d}''}^{X_{1},0}$.
\end{center}
\end{prop} 
{\em Proof.} We will first prove the statement 2. Consider the following diagram
\[\xymatrix{&K_{i-1,l-1}\ar[ld]_{\varphi_{\underline{d}',\underline{d}}}\\
K_{il}\ar[r]^{\varphi_{\underline{d},\underline{d}''}}&K_{i-1,l}\ar[u]_{\varphi_{\underline{d}'',\underline{d}'}}\ar[ld]^{\varphi_{\underline{d}'',\underline{d}'''}}\\
K_{i,l+1}\ar[u]^{\varphi_{\underline{d}''',\underline{d}}}&}\]
Notice that, by Proposition \ref{ultrarefinado es exacto}, item 2, elements $v_1,\ldots,v_{\beta}\in K_{il}$ exist satisfying that $\varphi_{\underline{d},\underline{d}''}(v_1)=u_1,\ldots,\varphi_{\underline{d},\underline{d}''}(v_{\beta})=u_{\beta}$. Since all sections of $\varphi_{\underline{d}',\underline{d}}(V_{\underline{d}'})\subseteq H^{0}(\mathcal{L}_{\underline{d}})$ vanish on $X_2$,
\begin{center}
$\varphi_{\underline{d}',\underline{d}}(V_{\underline{d}'})\cap \varphi_{\underline{d}''',\underline{d}}(V_{\underline{d}'''})\subseteq \varphi_{\underline{d}''',\underline{d}}(V_{\underline{d}'''})^{X_2,0}=\varphi_{\underline{d}''',\underline{d}}(V_{\underline{d}'''}^{X_2,0})=\varphi_{\underline{d}''',\underline{d}}(\varphi_{\underline{d}'',\underline{d}'''}(V_{\underline{d}''}))=\varphi_{\underline{d}'',\underline{d}}(V_{\underline{d}''})$,
\end{center}
where in the first equality we used Remark \ref{anulamiento} and in the second equality we used that $\varphi_{\underline{d}'',\underline{d}'''}(V_{\underline{d}''})=V_{\underline{d}'''}^{X_{2},0}$. On the other hand, $\varphi_{\underline{d}'',\underline{d}}(V_{\underline{d}''})=\varphi_{\underline{d}',\underline{d}}(\varphi_{\underline{d}'',\underline{d}'}(V_{\underline{d}''}))\subseteq \varphi_{\underline{d}',\underline{d}}(V_{\underline{d}'})$, as $\varphi_{\underline{d}'',\underline{d}'}(V_{\underline{d}''})\subseteq V_{\underline{d}'}$. Analogously, $\varphi_{\underline{d}'',\underline{d}}(V_{\underline{d}''})=\varphi_{\underline{d}''',\underline{d}}(\varphi_{\underline{d}'',\underline{d}'''}(V_{\underline{d}''}))\subseteq \varphi_{\underline{d}''',\underline{d}}(V_{\underline{d}'''})$, as $\varphi_{\underline{d}'',\underline{d}'''}(V_{\underline{d}''})\subseteq V_{\underline{d}'''}$. It follows that $\varphi_{\underline{d}'',\underline{d}}(V_{\underline{d}''})\subseteq \varphi_{\underline{d}',\underline{d}}(V_{\underline{d}'})\cap \varphi_{\underline{d}''',\underline{d}}(V_{\underline{d}'''})$, and hence $\varphi_{\underline{d}'',\underline{d}}(V_{\underline{d}''})=\varphi_{\underline{d}',\underline{d}}(V_{\underline{d}'})\cap \varphi_{\underline{d}''',\underline{d}}(V_{\underline{d}'''})$. Then
\begin{align}\label{exact ext1}
\text{dim}(\varphi_{\underline{d}',\underline{d}}(V_{\underline{d}'})+\varphi_{\underline{d}''',\underline{d}}(V_{\underline{d}'''}))=&\text{dim}\,\varphi_{\underline{d}',\underline{d}}(V_{\underline{d}'})+\text{dim}\,\varphi_{\underline{d}''',\underline{d}}(V_{\underline{d}'''})-\text{dim}\,\varphi_{\underline{d}'',\underline{d}}(V_{\underline{d}''}) \nonumber \\
=&(r+1-\text{dim}\,V_{\underline{d}'}^{X_{2}^{c},0})+(r+1-\text{dim}\,V_{\underline{d}'''}^{X_{3}^{c},0}) \nonumber \\
&-(r+1-\text{dim}\,V_{\underline{d}''}^{X_1,0}) \nonumber \\
=& r+1-(\text{dim}\,V_{\underline{d}'}^{X_{2}^{c},0}-(\text{dim}\,V_{\underline{d}''}^{X_1,0}-\text{dim}\,V_{\underline{d}'''}^{X_{3}^{c},0})). 
\end{align}
On the other hand, since $\varphi_{\underline{d}'',\underline{d}'''}(V_{\underline{d}''})=V_{\underline{d}'''}^{X_{2},0}$, we have
\begin{center}
$\varphi_{\underline{d}'',\underline{d}'''}(V_{\underline{d}''}^{X_1,0})=(\varphi_{\underline{d}'',\underline{d}'''}(V_{\underline{d}''}))^{X_1,0}=(V_{\underline{d}'''}^{X_{2},0})^{X_1,0}=V_{\underline{d}'''}^{X_{3}^{c},0}$,
\end{center}
where in the first equality we used Remark \ref{anulamiento}. Then 
\begin{equation}\label{exact ext2}
\text{dim}\,V_{\underline{d}''}^{X_1,0}-\text{dim}\,V_{\underline{d}'''}^{X_{3}^{c},0}=\text{dim}\,V_{\underline{d}''}^{X_{2}^{c},0}.
\end{equation}
Also, since $\varphi_{\underline{d}'',\underline{d}'}(V_{\underline{d}''})=V_{\underline{d}'}^{X_{3},0}$, we have
\begin{center}
$\varphi_{\underline{d}'',\underline{d}'}(V_{\underline{d}''}^{X_1,0})=(\varphi_{\underline{d}'',\underline{d}'}(V_{\underline{d}''}))^{X_1,0}=(V_{\underline{d}'}^{X_{3},0})^{X_1,0}=V_{\underline{d}'}^{X_{2}^{c},0}$,
\end{center}
and hence 
\begin{equation}\label{exact ext3}
\text{dim}\,V_{\underline{d}''}^{X_1,0}-\text{dim}\,V_{\underline{d}''}^{X_{3}^{c},0}=\text{dim}\,V_{\underline{d}'}^{X_{2}^{c},0}.
\end{equation}
It follows from (\ref{exact ext1}), (\ref{exact ext2}) and (\ref{exact ext3}) that 
\begin{align}\label{exact ext4}
\text{dim}(\varphi_{\underline{d}',\underline{d}}(V_{\underline{d}'})+\varphi_{\underline{d}''',\underline{d}}(V_{\underline{d}'''}))=& r+1-(\text{dim}\,V_{\underline{d}'}^{X_{2}^{c},0}-(\text{dim}\,V_{\underline{d}''}^{X_1,0}-\text{dim}\,V_{\underline{d}'''}^{X_{3}^{c},0})) \nonumber \\
=&r+1-(\text{dim}\,V_{\underline{d}'}^{X_{2}^{c},0}-\text{dim}\,V_{\underline{d}''}^{X_{2}^{c},0}) \nonumber \\
=&r+1-(\text{dim}\,V_{\underline{d}''}^{X_1,0}-\text{dim}\,V_{\underline{d}''}^{X_{3}^{c},0}-\text{dim}\,V_{\underline{d}''}^{X_{2}^{c},0}) \nonumber \\
=&r+1-\beta.
\end{align}
Then, to prove that $\text{dim}\,V_{\underline{d}}=r+1$, it suffices to show that 
\begin{center}
$(\varphi_{\underline{d}',\underline{d}}(V_{\underline{d}'})+\varphi_{\underline{d}''',\underline{d}}(V_{\underline{d}'''}))\cap \left\langle v_1,\ldots,v_{\beta}\right\rangle=0$.
\end{center}
Notice that
\begin{align}\label{exact ext5} 
\varphi_{\underline{d},\underline{d}''}(\varphi_{\underline{d}',\underline{d}}(V_{\underline{d}'})+\varphi_{\underline{d}''',\underline{d}}(V_{\underline{d}'''}))=&\varphi_{\underline{d},\underline{d}''}(\varphi_{\underline{d}',\underline{d}}(V_{\underline{d}'}))+\varphi_{\underline{d},\underline{d}''}(\varphi_{\underline{d}''',\underline{d}}(V_{\underline{d}'''})) \nonumber \\
=&\varphi_{\underline{d}',\underline{d}''}(V_{\underline{d}'})+\varphi_{\underline{d}''',\underline{d}''}(V_{\underline{d}'''}) \nonumber \\
=&V_{\underline{d}''}^{X_{3}^{c},0}+V_{\underline{d}''}^{X_{2}^{c},0} 
\end{align}
and
\begin{equation}\label{exact ext6}
\varphi_{\underline{d},\underline{d}''}(\left\langle v_1,\ldots,v_{\beta}\right\rangle)=\left\langle u_1,\ldots,u_{\beta}\right\rangle.
\end{equation}
Then 
\begin{equation*}
\varphi_{\underline{d},\underline{d}''}((\varphi_{\underline{d}',\underline{d}}(V_{\underline{d}'})+\varphi_{\underline{d}''',\underline{d}}(V_{\underline{d}'''}))\cap \left\langle v_1,\ldots,v_{\beta}\right\rangle)\subseteq (V_{\underline{d}''}^{X_{3}^{c},0}+V_{\underline{d}''}^{X_{2}^{c},0})\cap \left\langle u_1,\ldots,u_{\beta}\right\rangle=0,
\end{equation*}
where in the last equality we used that $V_{\underline{d}''}^{X_{1},0}=(V_{\underline{d}''}^{X_{2}^{c},0}\oplus V_{\underline{d}''}^{X_{3}^{c},0})\oplus \left\langle u_1,\ldots,u_{\beta}\right\rangle$. Therefore $\varphi_{\underline{d},\underline{d}''}((\varphi_{\underline{d}',\underline{d}}(V_{\underline{d}'})+\varphi_{\underline{d}''',\underline{d}}(V_{\underline{d}'''}))\cap \left\langle v_1,\ldots,v_{\beta}\right\rangle)=0$. On the other hand, as $u_1,\ldots,u_{\beta}$ are linearly independent and $\varphi_{\underline{d},\underline{d}''}(v_1)=u_1,\ldots,\varphi_{\underline{d},\underline{d}''}(v_{\beta})=u_{\beta}$, it follows that 
\begin{center}
$\varphi_{\underline{d},\underline{d}''}\big|_{\left\langle v_1,\ldots,v_{\beta}\right\rangle}:\left\langle v_1,\ldots,v_{\beta}\right\rangle\rightarrow \left\langle u_1,\ldots,u_{\beta}\right\rangle$ 
\end{center}
is an isomorphism. So $(\varphi_{\underline{d}',\underline{d}}(V_{\underline{d}'})+\varphi_{\underline{d}''',\underline{d}}(V_{\underline{d}'''}))\cap \left\langle v_1,\ldots,v_{\beta}\right\rangle=0$, and hence 
\begin{center}
$V_{\underline{d}}=(\varphi_{\underline{d}',\underline{d}}(V_{\underline{d}'})+\varphi_{\underline{d}''',\underline{d}}(V_{\underline{d}'''}))\oplus \left\langle v_1,\ldots,v_{\beta}\right\rangle$
\end{center}
is $r+1$-dimensional. Since $V_{\underline{d}'}$ and $V_{\underline{d}'''}$ are subspaces of $K_{i-1,l-1}$ and $K_{i,l+1}$, respectively, we have $\varphi_{\underline{d}',\underline{d}}(V_{\underline{d}'})\subseteq K_{il}$ and $\varphi_{\underline{d}''',\underline{d}}(V_{\underline{d}'''})\subseteq K_{il}$. Thus, as $\left\langle v_1,\ldots,v_{\beta}\right\rangle\subseteq K_{il}$ as well, it follows that $V_{\underline{d}}\subseteq K_{il}$.

Now, it follows from (\ref{exact ext5}) and (\ref{exact ext6}) that
\begin{equation}\label{exact ext7}
\varphi_{\underline{d},\underline{d}''}(V_{\underline{d}})=(V_{\underline{d}''}^{X_{3}^{c},0}+V_{\underline{d}''}^{X_{2}^{c},0})+\left\langle u_1,\ldots,u_{\beta}\right\rangle=V_{\underline{d}''}^{X_{1},0}.
\end{equation}
On the other hand, since $\varphi_{\underline{d}'',\underline{d}}(V_{\underline{d}''})\subseteq \varphi_{\underline{d}',\underline{d}}(V_{\underline{d}'})\subseteq V_{\underline{d}}$, it follows that
\begin{center}
$\text{dim}\,V_{\underline{d}''}^{X_{1},0}+\text{dim}\,\varphi_{\underline{d}'',\underline{d}}(V_{\underline{d}''})=\text{dim}\,V_{\underline{d}''}=r+1$,
\end{center}
so 
\begin{equation*}
\text{dim}\,\varphi_{\underline{d}'',\underline{d}}(V_{\underline{d}''})=r+1-\text{dim}\,V_{\underline{d}''}^{X_{1},0}
=\text{dim}\,V_{\underline{d}}-\text{dim}\,V_{\underline{d}''}^{X_{1},0}
=\text{dim}\,V_{\underline{d}}^{X_{1}^{c},0},
\end{equation*}
where the last equality follows from (\ref{exact ext7}). Thus, as $\varphi_{\underline{d}'',\underline{d}}(V_{\underline{d}''})\subseteq V_{\underline{d}}^{X_{1}^{c},0}$, it follows that $\varphi_{\underline{d}'',\underline{d}}(V_{\underline{d}''})=V_{\underline{d}}^{X_{1}^{c},0}$.

Now, we will show that 
$\varphi_{\underline{d}',\underline{d}}(V_{\underline{d}'})=V_{\underline{d}}^{X_{2},0}$ and $\varphi_{\underline{d},\underline{d}'}(V_{\underline{d}})=V_{\underline{d}'}^{X_{2}^{c},0}$. We have 
\begin{center}
$\varphi_{\underline{d},\underline{d}'}(V_{\underline{d}})=\varphi_{\underline{d}'',\underline{d}'}(\varphi_{\underline{d},\underline{d}''}(V_{\underline{d}}))=\varphi_{\underline{d}'',\underline{d}'}(V_{\underline{d}''}^{X_{1},0})=(\varphi_{\underline{d}'',\underline{d}'}(V_{\underline{d}''}))^{X_{1},0}=(V_{\underline{d}'}^{X_3,0})^{X_{1},0}=V_{\underline{d}'}^{X_{2}^{c},0}$,
\end{center}
where the second equality follows from (\ref{exact ext7}), and in the fourth equality we used that $\varphi_{\underline{d}'',\underline{d}'}(V_{\underline{d}''})=V_{\underline{d}'}^{X_3,0}$. Since $\varphi_{\underline{d}',\underline{d}}(V_{\underline{d}'})\subseteq V_{\underline{d}}$, it follows that
\begin{center}
$\text{dim}\,V_{\underline{d}'}^{X_{2}^{c},0}+\text{dim}\,\varphi_{\underline{d}',\underline{d}}(V_{\underline{d}'})=\text{dim}\,V_{\underline{d}'}=r+1$,
\end{center}
and hence
\begin{equation*}
\text{dim}\,\varphi_{\underline{d}',\underline{d}}(V_{\underline{d}'})=r+1-\text{dim}\,V_{\underline{d}'}^{X_{2}^{c},0}
=\text{dim}\,V_{\underline{d}}-\text{dim}\,V_{\underline{d}'}^{X_{2}^{c},0}
=\text{dim}\,V_{\underline{d}}^{X_{2},0},
\end{equation*}
where in the last equality we used that $\varphi_{\underline{d},\underline{d}'}(V_{\underline{d}})=V_{\underline{d}'}^{X_{2}^{c},0}$. Since $\varphi_{\underline{d}',\underline{d}}(V_{\underline{d}'})\subseteq V_{\underline{d}}^{X_{2},0}$, we get $\varphi_{\underline{d}',\underline{d}}(V_{\underline{d}'})=V_{\underline{d}}^{X_{2},0}$. The proof of the equalities $\varphi_{\underline{d}''',\underline{d}}(V_{\underline{d}'''})=V_{\underline{d}}^{X_{3},0}$ and $\varphi_{\underline{d},\underline{d}'''}(V_{\underline{d}})=V_{\underline{d}'''}^{X_{3}^{c},0}$ is analogous to that of $\varphi_{\underline{d}',\underline{d}}(V_{\underline{d}'})=V_{\underline{d}}^{X_{2},0}$ and $\varphi_{\underline{d},\underline{d}'}(V_{\underline{d}})=V_{\underline{d}'}^{X_{2}^{c},0}$. This proves statement 2.

Now, we will prove the statement 1. Notice that (\ref{exact ext3}) holds, as $\varphi_{\underline{d}'',\underline{d}'}(V_{\underline{d}''})=V_{\underline{d}'}^{X_{3},0}$. Then
\begin{align}\label{exact ext8}
\text{dim}\,\varphi_{\underline{d}',\underline{d}}(V_{\underline{d}'})=&r+1-\text{dim}\,V_{\underline{d}'}^{X_{2}^{c},0} \nonumber \\
=&r+1-(\text{dim}\,V_{\underline{d}''}^{X_1,0}-\text{dim}\,V_{\underline{d}''}^{X_{3}^{c},0}).
\end{align}
Since $i+l=d$, we have $\underline{d}''=(i-1,1,l)$, and hence $V_{\underline{d}''}^{X_{2}^{c},0}=0$. Thus 
\begin{equation}\label{exact ext9}
V_{\underline{d}''}^{X_{1},0}=(V_{\underline{d}''}^{X_{2}^{c},0}\oplus V_{\underline{d}''}^{X_{3}^{c},0})\oplus \left\langle u_1,\ldots,u_{\beta}\right\rangle=V_{\underline{d}''}^{X_{3}^{c},0}\oplus \left\langle u_1,\ldots,u_{\beta}\right\rangle.
\end{equation}
It follows from (\ref{exact ext8}) and (\ref{exact ext9}) that $\text{dim}\,\varphi_{\underline{d}',\underline{d}}(V_{\underline{d}'})=r+1-\beta$. To prove that $\text{dim}\,V_{\underline{d}}=r+1$, it suffices to show that 
\begin{center}
$\varphi_{\underline{d}',\underline{d}}(V_{\underline{d}'})\cap \left\langle v_1,\ldots,v_{\beta}\right\rangle=0$.
\end{center}
We have
\begin{align}\label{exact ext10} 
\varphi_{\underline{d},\underline{d}''}(\varphi_{\underline{d}',\underline{d}}(V_{\underline{d}'}))=\varphi_{\underline{d}',\underline{d}''}(V_{\underline{d}'})
=V_{\underline{d}''}^{X_{3}^{c},0}.
\end{align}
Then
\begin{equation*}
\varphi_{\underline{d},\underline{d}''}(\varphi_{\underline{d}',\underline{d}}(V_{\underline{d}'})\cap \left\langle v_1,\ldots,v_{\beta}\right\rangle)\subseteq V_{\underline{d}''}^{X_{3}^{c},0}\cap \left\langle u_1,\ldots,u_{\beta}\right\rangle=0,
\end{equation*}
where in the last equality we used (\ref{exact ext9}). Thus $\varphi_{\underline{d},\underline{d}''}(\varphi_{\underline{d}',\underline{d}}(V_{\underline{d}'})\cap \left\langle v_1,\ldots,v_{\beta}\right\rangle)=0$. Reasoning as in the proof of the statement 2, we get $\varphi_{\underline{d}',\underline{d}}(V_{\underline{d}'})\cap \left\langle v_1,\ldots,v_{\beta}\right\rangle=0$, and hence 
\begin{center}
$V_{\underline{d}}=\varphi_{\underline{d}',\underline{d}}(V_{\underline{d}'})\oplus \left\langle v_1,\ldots,v_{\beta}\right\rangle$
\end{center}
is $r+1$-dimensional. Reasoning as in the proof of the statement 2, we get $V_{\underline{d}}\subseteq K_{il}$.

Now, it follows from (\ref{exact ext9}) and (\ref{exact ext10}) that 
\begin{equation*}
\varphi_{\underline{d},\underline{d}''}(V_{\underline{d}})=V_{\underline{d}''}^{X_{3}^{c},0}+\left\langle u_1,\ldots,u_{\beta}\right\rangle=V_{\underline{d}''}^{X_{1},0}.
\end{equation*}
The proofs of the equalities $\varphi_{\underline{d}'',\underline{d}}(V_{\underline{d}''})=V_{\underline{d}}^{X_{1}^{c},0}$, $\varphi_{\underline{d}',\underline{d}}(V_{\underline{d}'})=V_{\underline{d}}^{X_{2},0}$ and $\varphi_{\underline{d},\underline{d}'}(V_{\underline{d}})=V_{\underline{d}'}^{X_{2}^{c},0}$ are the same as in the proof of the statement 2. So the statement 1 is shown.

Now, we will prove the statement 3. Notice that (\ref{exact ext2}) holds, as $\varphi_{\underline{d}'',\underline{d}'''}(V_{\underline{d}''})=V_{\underline{d}'''}^{X_{2},0}$. Therefore
\begin{align}\label{exact ext11}
\text{dim}\,\varphi_{\underline{d}''',\underline{d}}(V_{\underline{d}'''})=&r+1-\text{dim}\,V_{\underline{d}'''}^{X_{3}^{c},0} \nonumber \\
=&r+1-(\text{dim}\,V_{\underline{d}''}^{X_1,0}-\text{dim}\,V_{\underline{d}''}^{X_{2}^{c},0}).
\end{align}
Since $l=0$, we have $\underline{d}''=(i-1,d-i+1,0)$, and hence $V_{\underline{d}''}^{X_{3}^{c},0}=0$. Then
\begin{equation}\label{exact ext12}
V_{\underline{d}''}^{X_{1},0}=(V_{\underline{d}''}^{X_{2}^{c},0}\oplus V_{\underline{d}''}^{X_{3}^{c},0})\oplus \left\langle u_1,\ldots,u_{\beta}\right\rangle=V_{\underline{d}''}^{X_{2}^{c},0}\oplus \left\langle u_1,\ldots,u_{\beta}\right\rangle.
\end{equation}
It follows from (\ref{exact ext11}) and (\ref{exact ext12}) that $\text{dim}\,\varphi_{\underline{d}''',\underline{d}}(V_{\underline{d}'''})=r+1-\beta$. To prove that $\text{dim}\,V_{\underline{d}}=r+1$, it suffices to show that 
\begin{center}
$\varphi_{\underline{d}''',\underline{d}}(V_{\underline{d}'''})\cap \left\langle v_1,\ldots,v_{\beta}\right\rangle=0$.
\end{center}
We have
\begin{align}\label{exact ext13} 
\varphi_{\underline{d},\underline{d}''}(\varphi_{\underline{d}''',\underline{d}}(V_{\underline{d}'''}))=\varphi_{\underline{d}''',\underline{d}''}(V_{\underline{d}'''})
=V_{\underline{d}''}^{X_{2}^{c},0}.
\end{align}
Then
\begin{equation*}
\varphi_{\underline{d},\underline{d}''}(\varphi_{\underline{d}''',\underline{d}}(V_{\underline{d}'''})\cap \left\langle v_1,\ldots,v_{\beta}\right\rangle)\subseteq V_{\underline{d}''}^{X_{2}^{c},0}\cap \left\langle u_1,\ldots,u_{\beta}\right\rangle=0,
\end{equation*}
where in the last equality we used (\ref{exact ext12}). So $\varphi_{\underline{d},\underline{d}''}(\varphi_{\underline{d}''',\underline{d}}(V_{\underline{d}'''})\cap \left\langle v_1,\ldots,v_{\beta}\right\rangle)=0$, and reasoning as in the proof of the statement 2, we get $\varphi_{\underline{d}''',\underline{d}}(V_{\underline{d}'''})\cap \left\langle v_1,\ldots,v_{\beta}\right\rangle=0$, and hence 
\begin{center}
$V_{\underline{d}}=\varphi_{\underline{d}''',\underline{d}}(V_{\underline{d}'''})\oplus \left\langle v_1,\ldots,v_{\beta}\right\rangle$
\end{center}
is $r+1$-dimensional. Reasoning as in the proof of the statement 2, we get $V_{\underline{d}}\subseteq K_{il}$ as well.

Now, it follows from (\ref{exact ext12}) and (\ref{exact ext13}) that 
\begin{equation*}
\varphi_{\underline{d},\underline{d}''}(V_{\underline{d}})=V_{\underline{d}''}^{X_{2}^{c},0}+\left\langle u_1,\ldots,u_{\beta}\right\rangle=V_{\underline{d}''}^{X_{1},0}.
\end{equation*}
Reasoning as in the proof of the statement 2, we get the remaining equalities. So the statement 3 is shown, proving the proposition.    \hfill $\Box$

Now, we describe the method for the construction of exact extensions $\{(\mathcal{L}_{\underline{d}},V_{\underline{d}})\}_{\underline{d}}$. According to Remark \ref{ultrarefinado determina los V}, we necessarily have to do the construction in such a way that $V_{il}$ to be contained in $K_{il}$ for any $i\geq 0$ and $l\geq 0$ such that $i+l\leq d$. The idea is first to construct the subspaces $V_{il}$ for i=0, then for i=1, and so on, until $i=d-1$. For each $i\geq 1$, we first construct the subspaces $V_{il}$ for $l=d-i$, then for $l=d-i-1$, and so on, until $l=0$.

Suppose inductively that, for $i\geq 1$, the $r+1$-dimensional subspaces 
\begin{center}
$V_{i-1,l}\subseteq K_{i-1,l}$
\end{center}
have been constructed for $l=d-(i-1),\ldots,l=0$ in such a way that
\begin{center}
$\varphi_{\underline{d}',\underline{d}''}(V_{\underline{d}'})=V_{\underline{d}''}^{X_{3}^{c},0}$ and $\varphi_{\underline{d}'',\underline{d}'}(V_{\underline{d}''})=V_{\underline{d}'}^{X_{3},0}$ for $l=1,\ldots,d-(i-1)$,
\end{center}
where, as usual, $\underline{d}':=(i-1,d-(i-1)-l+1,l-1)$ and $\underline{d}'':=(i-1,d-(i-1)-l,l)$. We say that the subspaces $V_{i-1,l}$ satisfy the vertical exactness property.\\
Then we will construct $r+1$-dimensional subspaces $V_{il}\subseteq K_{il}$ for $l=d-i,\ldots,l=0$ in such a way that

(i) The subspaces $V_{il}$ satisfy the vertical exactness property.

(ii) 
\begin{center}
$\varphi_{\underline{d}',\underline{d}}(V_{\underline{d}'})=V_{\underline{d}}^{X_{2},0}$, $\varphi_{\underline{d},\underline{d}'}(V_{\underline{d}})=V_{\underline{d}'}^{X_{2}^{c},0}$ for $l=1,\ldots,d-i$, and\\
$\varphi_{\underline{d}'',\underline{d}}(V_{\underline{d}''})=V_{\underline{d}}^{X_{1}^{c},0}$, $\varphi_{\underline{d},\underline{d}''}(V_{\underline{d}})=V_{\underline{d}''}^{X_{1},0}$ for $l=0,\ldots,d-i$,
\end{center}
where $\underline{d}:=(i,d-i-l,l)$, $\underline{d}':=(i-1,d-(i-1)-l+1,l-1)$ and $\underline{d}'':=(i-1,d-(i-1)-l,l)$.\\
We inductively do the construction as follows:

{\em Step 1.} For $l=d-i$, the subspace $V_{il}$ is the subspace $V_{\underline{d}}$ defined in Proposition \ref{exact ext}, item 1.

{\em Step 2.} For $l=d-i-1,\ldots,l=1$, the subspace $V_{il}$ is the subspace $V_{\underline{d}}$ defined in Proposition \ref{exact ext}, item 2.

{\em Step 3.} For $l=0$, the subspace $V_{il}$ is the subspace $V_{\underline{d}}$ defined in Proposition \ref{exact ext}, item 3.

By Proposition \ref{exact ext}, the subspaces $V_{il}$ satisfy the properties (i) and (ii). Thus, it remains to construct the $r+1$-dimensional subspaces $V_{0l}\subseteq K_{0l}$ to satisfy the vertical exactness property, and verify that the exact limit linear series $\{(\mathcal{L}_{\underline{d}},V_{\underline{d}})\}_{\underline{d}}$ that we construct is in fact an extension.

By Proposition \ref{ultrarefinado es lls}, item 1, dim\,$K_{il}=r+1$ if $i\leq b_0$ or $l\leq b'_0$. Since $V_{X_1}, V_{X_2}$ and $V_{X_3}$ are linked, we have $V_{X_1}\subseteq K_{d0}$, $V_{X_2}\subseteq K_{00}$ and $V_{X_3}\subseteq K_{0d}$. It follows from dimension considerations, that $V_{X_1}=K_{d0}$, $V_{X_2}=K_{00}$ and $V_{X_3}=K_{0d}$. On the other hand, dim\,$K_{0l}=r+1$ if $0\leq l\leq d$. Thus, we define $V_{0l}:=K_{0l}$ for any nonnegative integer $l\leq d$. (Note that, for $l=0$ and $l=d$, the definition coincides with our fixed subspaces $V_{X_2}$ and $V_{X_3}$.)\\
Now, we will prove that 
\begin{center}
$\varphi_{\underline{d}''',\underline{d}}(V_{\underline{d}'''})=V_{\underline{d}}^{X_{3},0}$, $\varphi_{\underline{d},\underline{d}'''}(V_{\underline{d}})=V_{\underline{d}'''}^{X_{3}^{c},0}$ for $l=0,\ldots,d-1$,
\end{center}
where $\underline{d}:=(0,d-l,l)$ and $\underline{d}''':=(0,d-l-1,l+1)$. (Observe that this notation corresponds to the notation in Proposition \ref{exact ext} for $i=0$.)

Let $l$ be a nonnegative integer such that $l\leq d$. By Proposition \ref{ultrarefinado es exacto}, item 3, we have $\varphi_{\underline{d}''',\underline{d}}(K_{0,l+1})=K_{0l}^{X_{3},0}$. Then dim\,$K_{0,l+1}^{X_{3}^{c},0}+$dim\,$K_{0l}^{X_{3},0}=$dim\,$K_{0,l+1}$. On the other hand, since dim\,$K_{0,l+1}=r+1=$dim\,$K_{0,l}$, we have dim\,$K_{0,l+1}^{X_{3}^{c},0}+$dim\,$K_{0l}^{X_{3},0}=$dim\,$K_{0l}$. Now, dim\,$K_{0l}^{X_{3},0}+$dim\,$\varphi_{\underline{d},\underline{d}'''}(K_{0l})=$dim\,$K_{0l}$. Hence
\begin{center}
dim\,$\varphi_{\underline{d},\underline{d}'''}(K_{0l})=$dim\,$K_{0,l+1}^{X_{3}^{c},0}$,
\end{center}
and since $\varphi_{\underline{d},\underline{d}'''}(K_{0l})\subseteq K_{0,l+1}^{X_{3}^{c},0}$, we have $\varphi_{\underline{d},\underline{d}'''}(K_{0l})=K_{0,l+1}^{X_{3}^{c},0}$, proving that the subspaces $V_{0l}=K_{0l}$ satisfy the vertical exactness property.

Thus, we construct subspaces $V_{il}\subseteq K_{il}$ for $i=0,\ldots,d-1$ and $l=0,\ldots,d-i$. Now, since dim\,$K_{i0}=r+1$ for $i=0,\ldots,d-1$, we have that, by dimension considerations, $V_{i0}=K_{i0}$ for $i=0,\ldots,d-1$. On the other hand, the subspaces $\{K_{i0}\}_{i=0,\ldots,d}$ satisfy the following exactness property analogous to that of the subspaces $K_{0l}$
\begin{center}
$\varphi_{\underline{d}'',\underline{d}}(K_{i-1,0})=K_{i0}^{X_{1}^{c},0}$ and $\varphi_{\underline{d},\underline{d}''}(K_{i0})=K_{i-1,0}^{X_{1},0}$ for $i=1,\ldots,d$,
\end{center}
where $\underline{d}:=(i,d-i,0)$ and $\underline{d}'':=(i-1,d-i+1,0)$.\\
Thus, since $V_{X_1}=K_{d0}$, we get an exact limit linear series $\{(\mathcal{L}_{\underline{d}},V_{\underline{d}})\}_{\underline{d}}$ which is an extension of $\mathfrak{h}$.

Now, let $\{(\mathcal{L}_{\underline{d}},V_{\underline{d}})\}_{\underline{d}}$ be any exact extension. We will prove that the subspaces $V_{\underline{d}}$ are constructed by our method. Since dim\,$K_{0l}=r+1$ if $0\leq l\leq d$, we have $V_{0l}=K_{0l}$ for each integer $l$ such that $0\leq l\leq d$. Let $0<i<d$, $0\leq l\leq d-i$, and $\underline{d}:=(i,d-i-l,l)$. We will show that $V_{\underline{d}}$ is constructed by our method if $l>0$ and $i+l\leq d-1$. (The proofs of the cases $i+l=d$, $l=0$ are analogous.) Keep the notation of multidegrees used in Proposition \ref{exact ext}. By Remark \ref{ultrarefinado determina los V}, we have that $V_{\underline{d}}$, $V_{\underline{d}'}$, $V_{\underline{d}''}$ and $V_{\underline{d}'''}$ are $r+1$-dimensional subspaces of $K_{il}$, $K_{i-1,l-1}$, $K_{i-1,l}$ and $K_{i,l+1}$, respectively. Set $\beta:=\text{dim}\,V_{\underline{d}''}^{X_{1},0}-\text{dim}\,(V_{\underline{d}''}^{X_{2}^{c},0}\oplus V_{\underline{d}''}^{X_{3}^{c},0})$. As
$V_{\underline{d}}\supseteq \varphi_{\underline{d}',\underline{d}}(V_{\underline{d}'})$ and $V_{\underline{d}}\supseteq \varphi_{\underline{d}''',\underline{d}}(V_{\underline{d}'''})$, it follows that $V_{\underline{d}}\supseteq \varphi_{\underline{d}',\underline{d}}(V_{\underline{d}'})+\varphi_{\underline{d}''',\underline{d}}(V_{\underline{d}'''})$. By the proof of the statement 2 of Proposition \ref{exact ext}, we have
$\text{dim}\,(\varphi_{\underline{d}',\underline{d}}(V_{\underline{d}'})+\varphi_{\underline{d}''',\underline{d}}(V_{\underline{d}'''}))=r+1-\beta$, so 
\begin{center}
$V_{\underline{d}}=(\varphi_{\underline{d}',\underline{d}}(V_{\underline{d}'})+\varphi_{\underline{d}''',\underline{d}}(V_{\underline{d}'''}))\oplus \left\langle v_1,\ldots,v_{\beta}\right\rangle$,
\end{center}
for some $v_1,\ldots,v_{\beta}\in K_{il}$ which are linearly independent. 

Now, let $u_{1}:=\varphi_{\underline{d},\underline{d}''}(v_{1}),\ldots, u_{\beta}:=\varphi_{\underline{d},\underline{d}''}(v_{\beta})$. By the proof of the statement 2 of Proposition \ref{exact ext}, we have
\begin{center}
$\varphi_{\underline{d},\underline{d}''}(\varphi_{\underline{d}',\underline{d}}(V_{\underline{d}'})+\varphi_{\underline{d}''',\underline{d}}(V_{\underline{d}'''}))=V_{\underline{d}''}^{X_{3}^{c},0}+V_{\underline{d}''}^{X_{2}^{c},0}$, 
\end{center}
and hence
$V_{\underline{d}''}^{X_{1},0}=\varphi_{\underline{d},\underline{d}''}(V_{\underline{d}})=(V_{\underline{d}''}^{X_{3}^{c},0}+V_{\underline{d}''}^{X_{2}^{c},0})+\left\langle u_1,\ldots,u_{\beta}\right\rangle$. Since 
\begin{center}
$\beta=\text{dim}\,V_{\underline{d}''}^{X_{1},0}-\text{dim}\,(V_{\underline{d}''}^{X_{2}^{c},0}\oplus V_{\underline{d}''}^{X_{3}^{c},0})$, 
\end{center}
we necessarily have that $ u_1,\ldots,u_{\beta}$ are linearly independent and 
\begin{center}
$V_{\underline{d}''}^{X_{1},0}=(V_{\underline{d}''}^{X_{3}^{c},0}\oplus V_{\underline{d}''}^{X_{2}^{c},0})\oplus \left\langle u_1,\ldots,u_{\beta}\right\rangle$. 
\end{center}
This proves that our method constructs any exact extension.

\section{Unique exact extension}

In this section, we will show the conditions under which the exact extension is unique, and in this case, we will describe the scheme $\mathbb{P}(\mathfrak{g})$ for such a unique extension. Keeping the notation of multidegrees used in Proposition \ref{exact ext}, we have the following lemma.

\begin{lem}\label{dimension igual}
If $\varphi_{\underline{d}'',\underline{d}}(K_{i-1,l})=K_{il}^{X_{1}^{c},0}$, then {\em dim}\,$K_{il}=${\em dim}\,$K_{i-1,l}$.
\end{lem}
{\em Proof.} It follows from the hypothesis that 
\begin{center}
$\text{dim}\,K_{i-1,l}^{X_{1},0}+\text{dim}\,K_{il}^{X_{1}^{c},0}=\text{dim}\,K_{i-1,l}$.
\end{center}
On the other hand, by Proposition \ref{ultrarefinado es exacto}, item 2, we have $\varphi_{\underline{d},\underline{d}''}(K_{il})=K_{i-1,l}^{X_{1},0}$, and hence
\begin{center}
$\text{dim}\,K_{il}^{X_{1}^{c},0}+\text{dim}\,K_{i-1,l}^{X_{1},0}=\text{dim}\,K_{il}$.
\end{center}
Thus dim\,$K_{il}=$dim\,$K_{i-1,l}$.   \hfill $\Box$

\begin{lem}\label{dimension cae}
Let $i\geq 0$ and $l\geq 0$ such that $i+l\leq d$. Then, the following statements hold:

1. If $i>0$, then {\em dim}\,$K_{il}\geq$ {\em dim}\,$K_{i-1,l}$.

2. If $l>0$, then {\em dim}\,$K_{il}\geq$ {\em dim}\,$K_{i,l-1}$.
\end{lem}
{\em Proof.} We will only prove the statement 1, as the proof of the statement 2 is analogous. Let $\underline{d}:=(i,d-i-l,l)$ and $\underline{d}'':=(i-1,d-i-l+1,l)$. By Proposition \ref{ultrarefinado es exacto}, item 2, we have $\varphi_{\underline{d},\underline{d}''}(K_{il})=K_{i-1,l}^{X_{1},0}$. It follows that
\begin{center}
$\text{dim}\,K_{il}^{X_{1}^{c},0}+\text{dim}\,K_{i-1,l}^{X_{1},0}=\text{dim}\,K_{il}$.
\end{center}
On the other hand, since $\varphi_{\underline{d}'',\underline{d}}(K_{i-1,l})\subseteq K_{il}^{X_{1}^{c},0}$, we have
\begin{center}
$\text{dim}\,K_{i-1,l}=\text{dim}\,K_{i-1,l}^{X_{1},0}+\text{dim}\,\varphi_{\underline{d}'',\underline{d}}(K_{i-1,l})\leq \text{dim}\,K_{i-1,l}^{X_{1},0}+\text{dim}\,K_{il}^{X_{1}^{c},0}$.
\end{center}
It follows that dim\,$K_{il}\geq$ dim\,$K_{i-1,l}$, proving the statement 1 of the lemma.    \hfill $\Box$

\begin{thm}\label{unique}
The following statements are equivalent:\\
1. $\mathfrak{h}$ has a unique exact extension.\\
2. {\em dim\,}$K_{il}=r+1$ if $i+l\leq d$, $b_{j-1}<i\leq b_{j}$, $b'_{k-1}<l\leq b'_{k}$ and $j+k\leq r+1$.\\
3. $\mathfrak{h}$ has a unique extension.
\end{thm} 
{\em Proof.} First, for a fixed integer $j$ such that $1\leq j\leq r$, let us consider the following statement:\\
4. dim\,$K_{il}=r+1$ if $i+l\leq d$, $b_{j-1}<i\leq b_{j}$, $b'_{k-1}<l\leq b'_{k}$ and $j+k\leq r+1$.\\ 
We will prove that, for that fixed integer $j$, the statement 4 implies the following statement:\\
5. Let $e_{j0},\ldots,e_{j,r-j}$ be the orders of vanishing of $V_{X_2}(-b_{j}A)$ at $B$. Then $e_{j,r-j}=b'_{r-j}$.

Indeed, assume statement 4 holds. by the proof of the first four cases of Proposition \ref{ultrarefinado es lls}, we have that, for $b_{j-1}<i\leq b_{j}$ and $b'_{k-1}<l\leq b'_{k}$ such that $i+l\leq d$, $\text{dim}\,K_{il}=r+1$ if and only if $\text{dim}\,V_{X_2}(-iA-lB)=r+1-j-k$. On the other hand, by Remark \ref{sumaordenes}, $b_{j}+b'_{r-j}\leq d$. Thus $\text{dim}\,K_{b_{j},b'_{r-j}}=r+1$, and it follows that 
\begin{center}
$\text{dim}\,V_{X_2}(-b_{j}A-b'_{r-j}B)=r+1-j-(r-j)=1$.
\end{center}
If $b_{j}+b'_{r-j}+1\leq d$, then, by hypothesis, $\text{dim}\,K_{b_{j},b'_{r-j}+1}=r+1$, as $b'_{r-j}<b'_{r-j}+1\leq b'_{r+1-j}$. So $\text{dim}\,V_{X_2}(-b_{j}A-(b'_{r-j}+1)B)=r+1-j-(r+1-j)=0$, and it follows that $b'_{r-j}=e_{r-j}$. Now, if $b_{j}+b'_{r-j}+1>d$, then $V_{X_2}(-b_{j}A-(b'_{r-j}+1)B)=0$ as well, and hence $b'_{r-j}=e_{j,r-j}$. So statement 5 holds. 

Now, we will prove that the statement 1 implies the statement 2. Assume statement 1 holds. Let $\{(\mathcal{L}_{\underline{d}},V_{\underline{d}})\}_{\underline{d}}$ be the unique exact extension. We claim that, for $0<i<d$ and $0\leq l\leq d-i$, 
\begin{center}
$K_{il}^{X_{1}^{c},0}=\varphi_{\underline{d}'',\underline{d}}(V_{\underline{d}''})$ if $\beta>0$,
\end{center}
where $\underline{d}:=(i,d-i-l,l)$, $\underline{d}'':=(i-1,d-i-l+1,l)$ and $\beta:=$dim\,$V_{\underline{d}''}^{X_{1},0}-$dim$(V_{\underline{d}''}^{X_{2}^{c},0}\oplus V_{\underline{d}''}^{X_{3}^{c},0})$.

Indeed, we will prove the claim for $l>0$ and $i+l\leq d-1$, as the remaining cases are analogous. Assume $\beta>0$. In Section \ref{extension} we saw that 
\begin{center}
$V_{\underline{d}}=(\varphi_{\underline{d}',\underline{d}}(V_{\underline{d}'})+\varphi_{\underline{d}''',\underline{d}}(V_{\underline{d}'''}))\oplus \left\langle v_1,\ldots,v_{\beta}\right\rangle$,
\end{center}
where $\underline{d}':=(i-1,d-i-l+2,l-1)$, $\underline{d}''':=(i,d-i-l-1,l+1)$ and $v_1,\ldots,v_{\beta}\in K_{il}$ satisfy that 
$V_{\underline{d}''}^{X_{1},0}=(V_{\underline{d}''}^{X_{2}^{c},0}\oplus V_{\underline{d}''}^{X_{3}^{c},0})\oplus \left\langle u_1,\ldots,u_{\beta}\right\rangle$, where $u_{1}:=\varphi_{\underline{d},\underline{d}''}(v_1),\ldots,u_{\beta}:=\varphi_{\underline{d},\underline{d}''}(v_{\beta})$. Suppose that $K_{il}^{X_{1}^{c},0}$ is not contained in $V_{\underline{d}}$. Let $\widetilde{v}\in K_{il}^{X_{1}^{c},0}\setminus V_{\underline{d}}$, and set 
\begin{center}
$\widetilde{V}_{\underline{d}}:=(\varphi_{\underline{d}',\underline{d}}(V_{\underline{d}'})+\varphi_{\underline{d}''',\underline{d}}(V_{\underline{d}'''}))\oplus \left\langle v_1+\widetilde{v},\ldots,v_{\beta}\right\rangle$.
\end{center}
We have that $\widetilde{V}_{\underline{d}}\subseteq K_{il}$ and $\varphi_{\underline{d},\underline{d}''}(v_1+\widetilde{v})=u_{1},\ldots,\varphi_{\underline{d},\underline{d}''}(v_{\beta})=u_{\beta}$. Now, as $v_1\in V_{\underline{d}}$ and $\widetilde{v}\notin V_{\underline{d}}$, it follows that $v_1+\widetilde{v}\notin V_{\underline{d}}$, and hence $\widetilde{V}_{\underline{d}}\neq V_{\underline{d}}$. However, by the method of Section \ref{extension}, this allows us to construct an exact extension which is different from the unique exact extension, a contradiction. Thus $K_{il}^{X_{1}^{c},0}\subseteq V_{\underline{d}}$, and hence $K_{il}^{X_{1}^{c},0}=V_{\underline{d}}^{X_{1}^{c},0}=\varphi_{\underline{d}'',\underline{d}}(V_{\underline{d}''})$. So our claim is established.

Now, we will prove the statement 2 by induction on $j$. Let $l_0$ be the largest order of vanishing of $V_{X_2}(-b_{1}A)$ at $B$. Notice that $b_{1}+l_{0}\leq d$, as $V_{X_2}(-b_{1}A-l_{0}B)\neq 0$. Also, we have $l_{0}\leq b'_{r}$, as $l_{0}$ is necessarily an order of vanishing of $V_{X_2}$ at $B$. By definition of $l_{0}$, and since $V_{X_2}(-(b_{0}+1)A)=V_{X_2}(-b_{1}A)$, we get
\begin{center}
dim\,$V_{X_2}(-(b_{0}+1)A-l_{0}B)=1$ and dim\,$V_{X_2}(-(b_{0}+1)A-(l_{0}+1)B)=0$. 
\end{center}
Thus, since $K_{il}^{X_{2}^{c},0}\cong V_{X_2}(-(i+1)A-(l+1)B)$ for any $i\geq 0$ and $l\geq 0$ such that $i+l\leq d$, we get 
\begin{align}\label{unique 1}
\text{dim}\,K_{b_{0},l_{0}-1}^{X_{2}^{c},0}-\text{dim}\,K_{b_{0},l_{0}}^{X_{2}^{c},0}=1-0=1\,\,\text{if}\,\,l_{0}>0.
\end{align}
Now, set $i:=b_{0}+1$ and $l:=l_{0}$, keep the notation of multidegrees used in Proposition \ref{exact ext} and set $\beta:=$dim\,$V_{\underline{d}''}^{X_{1},0}-$dim$(V_{\underline{d}''}^{X_{2}^{c},0}\oplus V_{\underline{d}''}^{X_{3}^{c},0})$. We will prove that $\beta>0$. Suppose first that $l>0$. It follows from the proofs of the statements 1 and 2 of Proposition \ref{exact ext} that $\beta=$dim\,$V_{\underline{d}'}^{X_{2}^{c},0}-$dim\,$V_{\underline{d}''}^{X_{2}^{c},0}$. On the other hand, dim\,$K_{b_{0},l_{0}}=r+1$ and dim\,$K_{b_{0},l_{0}-1}=r+1$, so $V_{\underline{d}''}=K_{b_{0},l_{0}}$ and $V_{\underline{d}'}=K_{b_{0},l_{0}-1}$. Therefore, by (\ref{unique 1}),
\begin{center}
$\beta=$dim\,$V_{\underline{d}'}^{X_{2}^{c},0}-$dim\,$V_{\underline{d}''}^{X_{2}^{c},0}$=dim\,$K_{b_{0},l_{0}-1}^{X_{2}^{c},0}-$dim\,$K_{b_{0},l_{0}}^{X_{2}^{c},0}=1$,
\end{center}
so $\beta>0$. Suppose now that $l=0$. It follows from the proof of the statement 3 of Proposition \ref{exact ext} that $\beta=$dim\,$V_{\underline{d}''}^{X_{1},0}-$dim\,$V_{\underline{d}''}^{X_{2}^{c},0}$. On the other hand, as we saw, $V_{\underline{d}''}=K_{b_{0},l_{0}}$. Since
$\text{dim}\,K_{b_{0},l_{0}}^{X_{2}^{c},0}=0$ and
\begin{equation*}
\begin{split}
\text{dim}\,K_{b_{0},l_{0}}^{X_1,0}&=\text{dim}\,V_{X_2}(-(b_{0}+1)A-l_{0}B)+\text{dim}\,V_{X_3}(-(d-l_{0}+1)B) \\
&=1+\text{dim}\,V_{X_3}(-(d-0+1)B)=1,
\end{split}
\end{equation*}
we get 
\begin{center}
$\beta=$dim\,$V_{\underline{d}''}^{X_{1},0}-$dim\,$V_{\underline{d}''}^{X_{2}^{c},0}=\text{dim}\,K_{b_{0},l_{0}}^{X_1,0}-\text{dim}\,K_{b_{0},l_{0}}^{X_{2}^{c},0}=1-0=1$.
\end{center}
Thus, in any case, $\beta>0$. It follows from the claim that $K_{il}^{X_{1}^{c},0}=\varphi_{\underline{d}'',\underline{d}}(V_{\underline{d}''})\subseteq \varphi_{\underline{d}'',\underline{d}}(K_{i-1,l})$, and hence $K_{il}^{X_{1}^{c},0}=\varphi_{\underline{d}'',\underline{d}}(K_{i-1,l})$. Then, by lemma \ref{dimension igual}, we get 
\begin{align}\label{unique 2}
\text{dim}\,K_{il}=\text{dim}\,K_{i-1,l}=\text{dim}\,K_{b_{0},l_{0}}=r+1.
\end{align}
On the other hand, notice that, for $\tilde{l}>l_{0}$, $V_{X_2}(-b_{0}A-\tilde{l}B)=V_{X_2}(-(b_{0}+1)A-\tilde{l}B)$ if and only if $V_{X_2}(-\tilde{l}B)=V_{X_2}(-(b_{0}+1)A-\tilde{l}B)$, i.e., if and only if $V_{X_2}(-\tilde{l}B)=0$, that is, $\tilde{l}>b'_{r}$. Thus, $b_{0}$ is an order of vanishing of $V_{X_2}(-\tilde{l}B)$ if $l_{0}<\tilde{l}\leq b'_{r}$. Then, by Proposition \ref{exacto otro sentido}, item 1, if $l_{0}<\tilde{l}\leq b'_{r}$ and $i+\tilde{l}\leq d$, $K_{i\tilde{l}}^{X_{1}^{c},0}=\varphi_{\widetilde{\underline{d}}'',\widetilde{\underline{d}}}(K_{i-1,\tilde{l}})$, where $\widetilde{\underline{d}}:=(i,d-i-\tilde{l},\tilde{l})$ and $\widetilde{\underline{d}}'':=(i-1,d-i-\tilde{l}+1,\tilde{l})$. It follows from lemma \ref{dimension igual} that 
\begin{align}\label{unique 3}
\text{dim}\,K_{i\tilde{l}}=\text{dim}\,K_{i-1,\tilde{l}}=\text{dim}\,K_{b_{0},\tilde{l}}=r+1\, \text{if}\,\, l_{0}<\tilde{l}\leq b'_{r}\,\, \text{and}\,\, i+\tilde{l}\leq d.
\end{align}
Therefore, by (\ref{unique 2}), (\ref{unique 3}), lemma \ref{dimension cae}, item 2 and Proposition \ref{ultrarefinado es lls}, item 1, 
\begin{align}\label{unique 4}
\text{dim}\,K_{i\tilde{l}}=r+1 \,\text{if}\,\, 0\leq \tilde{l}\leq b'_{r}\, \text{and}\,\, i+\tilde{l}\leq d.
\end{align}
On the other hand, by Proposition \ref{exacto otro sentido}, item 1, if $b_{0}+1<\tilde{i}\leq b_{1}$, $\tilde{l}\geq 0$ and $\tilde{i}+\tilde{l}\leq d$, $K_{\tilde{i}\tilde{l}}^{X_{1}^{c},0}=\varphi_{\widetilde{\underline{d}}'',\widetilde{\underline{d}}}(K_{\tilde{i}-1,\tilde{l}})$, where $\widetilde{\underline{d}}:=(\tilde{i},d-\tilde{i}-\tilde{l},\tilde{l})$ and $\widetilde{\underline{d}}'':=(\tilde{i}-1,d-\tilde{i}-\tilde{l}+1,\tilde{l})$. Then, it follows from lemma \ref{dimension igual} that 
\begin{align}\label{unique 5}
\text{dim}\,K_{\tilde{i}\tilde{l}}=\text{dim}\,K_{\tilde{i}-1,\tilde{l}}=\ldots=\text{dim}\,K_{b_{0}+1,\tilde{l}}=r+1,
\end{align}
$\text{if}\,\,b_{0}+1<\tilde{i}\leq b_{1},\, 0\leq\tilde{l}\leq b'_{r}\,\,\text{and}\,\, \tilde{i}+\tilde{l}\leq d$. Thus, (\ref{unique 4}) and (\ref{unique 5}) prove the case $j=1$. Now, suppose by induction that, for a certain $2\leq j<r$, 
\begin{center}
dim\,$K_{\tilde{i}\tilde{l}}=r+1$ if $\tilde{i}+\tilde{l}\leq d$, $b_{j-1}<\tilde{i}\leq b_{j}$, $b'_{k-1}<\tilde{l}\leq b'_{k}$ and $j+k\leq r+1$.
\end{center}
Then, since the statement 4 implies the statement 5, $l_{j-1}:=b'_{r-j}$ is the largest order of vanishing of $V_{X_2}(-b_{j}A)$ at $B$. Let $l_{j}$ be the largest order of vanishing of $V_{X_2}(-b_{j+1}A)$ at $B$. Notice that $b_{j+1}+l_{j}\leq d$, as $V_{X_2}(-b_{j+1}A-l_{j}B)\neq 0$. Also, we have $l_{j}\leq l_{j-1}$, i.e., $l_{j}\leq b'_{r-j}$. By definition of $l_{j}$, and since $V_{X_2}(-(b_{j}+1)A)=V_{X_2}(-b_{j+1}A)$, we get
\begin{center}
dim\,$V_{X_2}(-(b_{j}+1)A-l_{j}B)=1$ and dim\,$V_{X_2}(-(b_{j}+1)A-(l_{j}+1)B)=0$,
\end{center}
and hence
\begin{align}\label{unique 6}
\text{dim}\,K_{b_{j},l_{j}-1}^{X_{2}^{c},0}-\text{dim}\,K_{b_{j},l_{j}}^{X_{2}^{c},0}=1-0=1\,\,\text{if}\,\,l_{j}>0.
\end{align}
Now, we proceed as in the first case. Set $i:=b_{j}+1$ and $l:=l_{j}$, keep the notation of multidegrees used in Proposition \ref{exact ext} and set $\beta:=$dim\,$V_{\underline{d}''}^{X_{1},0}-$dim$(V_{\underline{d}''}^{X_{2}^{c},0}\oplus V_{\underline{d}''}^{X_{3}^{c},0})$. We will show that $\beta>0$. Suppose first that $l>0$. As we saw, $\beta=$dim\,$V_{\underline{d}'}^{X_{2}^{c},0}-$dim\,$V_{\underline{d}''}^{X_{2}^{c},0}$. On the other hand, by induction, dim\,$K_{b_{j},l_{j}}=r+1$ and dim\,$K_{b_{j},l_{j}-1}=r+1$, as $b_{j}+l_{j}\leq d$ and $l_{j}\leq b'_{r-j}$, so $V_{\underline{d}''}=K_{b_{j},l_{j}}$ and $V_{\underline{d}'}=K_{b_{j},l_{j}-1}$. Therefore, by (\ref{unique 6}),
\begin{center}
$\beta=$dim\,$V_{\underline{d}'}^{X_{2}^{c},0}-$dim\,$V_{\underline{d}''}^{X_{2}^{c},0}$=dim\,$K_{b_{j},l_{j}-1}^{X_{2}^{c},0}-$dim\,$K_{b_{j},l_{j}}^{X_{2}^{c},0}=1$,
\end{center}
so $\beta>0$. Suppose now that $l=0$. We have $\beta=$dim\,$V_{\underline{d}''}^{X_{1},0}-$dim\,$V_{\underline{d}''}^{X_{2}^{c},0}$ and $V_{\underline{d}''}=K_{b_{j},l_{j}}$. Since
$\text{dim}\,K_{b_{j},l_{j}}^{X_{2}^{c},0}=0$ and
\begin{equation*}
\begin{split}
\text{dim}\,K_{b_{j},l_{j}}^{X_1,0}&=\text{dim}\,V_{X_2}(-(b_{j}+1)A-l_{j}B)+\text{dim}\,V_{X_3}(-(d-l_{j}+1)B) \\
&=1+\text{dim}\,V_{X_3}(-(d-0+1)B)=1,
\end{split}
\end{equation*}
we get 
\begin{center}
$\beta=$dim\,$V_{\underline{d}''}^{X_{1},0}-$dim\,$V_{\underline{d}''}^{X_{2}^{c},0}=\text{dim}\,K_{b_{j},l_{j}}^{X_1,0}-\text{dim}\,K_{b_{j},l_{j}}^{X_{2}^{c},0}=1-0=1$.
\end{center}
Thus, in any case, $\beta>0$. Reasoning as in the case $j=1$, we get 
\begin{align}\label{unique 7}
\text{dim}\,K_{il}=\text{dim}\,K_{b_{j},l_{j}}=r+1.
\end{align}
On the other hand, notice that, for $\tilde{l}>l_{j}$, $V_{X_2}(-b_{j}A-\tilde{l}B)=V_{X_2}(-(b_{j}+1)A-\tilde{l}B)$ if and only if $V_{X_2}(-b_{j}A-\tilde{l}B)=0$, i.e., if and only if $\tilde{l}>l_{j-1}=b'_{r-j}$. Thus, $b_{j}$ is an order of vanishing of $V_{X_2}(-\tilde{l}B)$ if $l_{j}<\tilde{l}\leq b'_{r-j}$. Then, by Proposition \ref{exacto otro sentido}, item 1, if $l_{j}<\tilde{l}\leq b'_{r-j}$ and $i+\tilde{l}\leq d$, $K_{i\tilde{l}}^{X_{1}^{c},0}=\varphi_{\widetilde{\underline{d}}'',\widetilde{\underline{d}}}(K_{i-1,\tilde{l}})$, where $\widetilde{\underline{d}}:=(i,d-i-\tilde{l},\tilde{l})$ and $\widetilde{\underline{d}}'':=(i-1,d-i-\tilde{l}+1,\tilde{l})$. By induction and lemma \ref{dimension igual},
\begin{align}\label{unique 8}
\text{dim}\,K_{i\tilde{l}}=\text{dim}\,K_{i-1,\tilde{l}}=\text{dim}\,K_{b_{j},\tilde{l}}=r+1\, \text{if}\,\, l_{j}<\tilde{l}\leq b'_{r-j}\,\, \text{and}\,\, i+\tilde{l}\leq d.
\end{align}
By (\ref{unique 7}) and (\ref{unique 8}),
\begin{center}
$\text{dim}\,K_{i\tilde{l}}=r+1 \,\text{if}\,\, 0\leq \tilde{l}\leq b'_{r-j}\, \text{and}\,\, i+\tilde{l}\leq d$.
\end{center}
Reasoning as in the case $j=1$, we get
\begin{center}
$\text{dim}\,K_{\tilde{i}\tilde{l}}=r+1 \,\text{if}\,\, b_{j}<\tilde{i}\leq b_{j+1},\,\, 0\leq \tilde{l}\leq b'_{r-j}\, \text{and}\,\, \tilde{i}+\tilde{l}\leq d$.
\end{center}
This finishes the proof of the statement 2.

Now, we will prove that the statement 2 implies the statement 3. Assume statement 2 holds. Let $\{(\mathcal{L}_{\underline{d}},V_{\underline{d}})\}_{\underline{d}}$ be an extension. Then, since $b_{j}+b'_{r-j}\leq d$, we have 
\begin{center}
$V_{il}=K_{il}$ if $0<j\leq r$, $b_{j-1}<i\leq b_{j}$ and $0\leq l\leq b'_{r-j}$. 
\end{center}
Also, 
\begin{center}
$V_{il}=K_{il}$ if $i\leq b_{0}$ or $l\leq b'_{0}$.
\end{center}
On the other hand, since the statement 4 implies the statement 5, for each $0<j\leq r$, $b'_{r-j}$ is the largest order of vanishing of $V_{X_2}(-b_{j}A)=V_{X_2}(-(b_{j-1}+1)A)$ at $B$. Then $V_{X_2}(-(b_{j-1}+1)A-(b'_{r-j}+1)B)=0$, and hence $V_{X_2}(-(i+1)A-(l+1)B)=0$ if $i\geq b_{j-1}$ and $l\geq b'_{r-j}$. It follows that $K_{il}^{X_{2}^{c},0}=0$ if $i\geq b_{j-1}$, $l\geq b'_{r-j}$ and $i+l\leq d$. Thus, for $i>b_{j-1}$, $l>b'_{r-j}$ and $i+l\leq d$, $K_{i-1,l-1}^{X_{2}^{c},0}=0$, and hence $V_{i-1,l-1}^{X_{2}^{c},0}=0$, implying that 
\begin{center}
$\varphi_{(i-1,d-i-l+2,l-1),(i,d-i-l,l)}(V_{i-1,l-1})=V_{il}$ if $i>b_{j-1}$, $l>b'_{r-j}$ and $i+l\leq d$.
\end{center}
(In this case, $V_{il}=V_{il}^{X_{2},0}$.) It follows that the extension is unique, proving the statement 3.

Finally, by the method of Section \ref{extension}, there exists at least one exact extension, so the statement 3 implies the statement 1. This finishes the proof of the theorem.   \hfill $\Box$

\begin{rem}\label{vacio}
{\em Suppose that $\mathfrak{h}$ has a unique exact extension and let $\{(\mathcal{L}_{\underline{d}},V_{\underline{d}})\}_{\underline{d}}$ be its unique exact extension. Note that, in the proof of Theorem \ref{unique}, we saw that
$V_{\underline{d}}=V_{\underline{d}}^{X_{2},0}$ for any $\underline{d}=(i,d-i-l,l)$ with $i>b_{j-1}$, $l>b'_{r-j}$ and $i+l\leq d$.
}
\end{rem}
\begin{rem}\label{ordenes}
{\em Suppose that $\mathfrak{h}$ has a unique exact extension and let $\{(\mathcal{L}_{\underline{d}},V_{\underline{d}})\}_{\underline{d}}$ be its unique exact extension. It follows from the proof of Theorem \ref{unique} that, for $0\leq j\leq r$, $b'_0,\ldots,b'_{r-j}$ are the orders of vanishing of $V_{X_2}(-b_{j}A)$ at $B$. Analogously, for $0\leq k\leq r$, $b_0,\ldots,b_{r-k}$ are the orders of vanishing of $V_{X_2}(-b'_{k}B)$ at $A$.
}
\end{rem}

\textbf{Acknowledgements.} The author would like to thank Eduardo Esteves for several helpful discussions.

\bibliographystyle{alpha}

\end{document}